\documentclass{article}

\usepackage{
 amsmath,
 amsxtra,
 amsthm,
 amssymb,
 etex,
 mathrsfs,
 stmaryrd,
pdflscape
 }
\usepackage[multiple]{footmisc} 
\usepackage[all]{xy}

\newtheorem{theorem}{Theorem}[section]
\newtheorem{lemma}[theorem]{Lemma}

\newtheorem{proposition}[theorem]{Proposition}
\newtheorem{corollary}[theorem]{Corollary}

\theoremstyle{definition}
\newtheorem{defn}[theorem]{Definition}
\newtheorem{remark}[theorem]{Remark}

\newcommand{\bd}{\begin{defn}}
\newcommand{\ed}{\end{defn}}
\newcommand{\bl}{\begin{lemma}}
\newcommand{\el}{\end{lemma}}
\newcommand{\bp}{\begin{proposition}}
\newcommand{\ep}{\end{proposition}}
\newcommand{\bt}{\begin{theorem}}
\newcommand{\et}{\end{theorem}}
\newcommand{\bc}{\begin{corollary}}
\newcommand{\ec}{\end{corollary}}
\newcommand{\br}{\begin{remark}}
\newcommand{\er}{\end{remark}}
\newcommand{\ba}{\begin{array}}
\newcommand{\ea}{\end{array}}
\newcommand{\bpf}{\begin{proof}}
\newcommand{\epf}{\end{proof}}

\newcommand{\Z}{\mathbb{Z}}
\newcommand{\Q}{\mathbb{Q}}
\newcommand{\Zp}{\mathbb{Z}_{p}}
\newcommand{\Qp}{\mathbb{Q}_{p}}
\newcommand{\Fp}{\mathbb{F}_{p}}
\newcommand{\Op}{\mathcal{O}}

\newcommand{\al}{\alpha}
\newcommand{\be}{\beta}
\newcommand{\Ga}{\Gamma}
\newcommand{\ga}{\gamma}

\newcommand{\la}{\lambda}

\newcommand{\Si}{\Sigma}

\DeclareMathOperator{\Sel}{Sel} \DeclareMathOperator{\Gal}{Gal}
\DeclareMathOperator{\Hom}{Hom} \DeclareMathOperator{\rank}{rank}
\DeclareMathOperator{\corank}{corank}
\DeclareMathOperator{\Ext}{Ext} \DeclareMathOperator{\Tor}{Tor}

\newcommand{\cyc}{\mathrm{cyc}}
\newcommand{\ch}{\mathrm{char}}

\newcommand{\M}{\mathfrak{M}}
\newcommand{\N}{\mathfrak{N}}

\newcommand{\ot}{\otimes}
\newcommand{\ilim}{\displaystyle \mathop{\varinjlim}\limits}
\newcommand{\plim}{\displaystyle \mathop{\varprojlim}\limits}

\newcommand{\coker}{\mathrm{coker}\,}

\newcommand{\lra}{\longrightarrow}

\newcommand{\ps}[1]{\llbracket #1 \rrbracket}

  \DeclareFontFamily{U}{wncy}{}
  \DeclareFontShape{U}{wncy}{m}{n}{<->wncyr10}{}
  \DeclareSymbolFont{mcy}{U}{wncy}{m}{n}
  \DeclareMathSymbol{\sha}{\mathord}{mcy}{"58}

\begin{document}
\title{On characteristic elements modulo $p$ in non-commutative Iwasawa theory}
 \author{
  Meng Fai Lim\footnote{School of Mathematics and Statistics, Hubei Key Laboratory of Mathematical Sciences,
Central China Normal University, Wuhan, 430079, P.R.China.
 E-mail: \texttt{limmf@ccnu.edu.cn}} \quad
 Chao Qin\footnote{College of Mathematical Sciences,
Harbin Engineering University,
Harbin, 150001, P.R.China.
 E-mail: \texttt{qinchao@hrbeu.edu.cn}}\footnote{The author’s research is supported by the National Natural Science Foundation of China under Grant No. 12001546,
Heilongjiang Province under Grant No. 3236330122.}} 

\date{}
\maketitle

\begin{abstract} \footnotesize
\noindent
 Coates, Fukaya, Kato, Sujatha and Venjakob come up with a procedure of attaching suitable
 characteristic element to Selmer groups defined over a non-commutative $p$-adic Lie extension, which is subsequently refined by Burns and Venjakob. By their construction, these
 characteristic elements are realized as elements in an appropriate localized $K_1$-group. In this paper, we will introduce a notion of modulo $p$ for these elements and study some of their properties. As an application, we study the Greenberg Selmer group of a tensor product of modular forms, where $p$ is an Eisenstein prime for one of these forms. 


\end{abstract}

\section{Introduction}

Throughout the paper, $p$ will always denote an odd prime. Let $\Gamma$ denote the multiplicative group which is isomorphic to the additive group of the $p$-adic ring of integers. In Iwasawa theory, the completed group algebra $\Zp\ps{\Gamma}$ naturally comes into play in the formulation of the main conjecture (for instances, see \cite{Iw73, G89, K}), where one usually attaches to an Iwasawa theoretical object its characteristic element which lives in 
the ring $\Zp\ps{\Gamma}$. As in number theory, one has to work $\alpha$ mod $p$ for a given element $\alpha\in \Zp\ps{\Gamma}$. The latter can be thought as understanding $\pi_{\Gamma}(\alpha)$, where $\pi_{\Gamma}$ is the natural projection $\Zp\ps{\Ga}\twoheadrightarrow \Fp\ps{\Ga}$, where $\Fp$ is the finite field with cardinality $p$. In particular, for instances, if $M$ is a $\Zp\ps{\Gamma}$-module which is also finitely generated over $\Zp$, and $f(x)$ is the characteristic polynomial of $M$ in the sense of \cite[Definition 5.3.9]{NSW}, we then have $\pi_{\Gamma}(f(x)) = T^{\rank_{\Zp}(M)}$. This viewpoint has been taken up, most notably, in the works of Greenberg-Vatsal \cite{GV} and Emerton-Pollack-Weston \cite{EPW} which has since inspired many later works (see \cite{AAS, Ha, JMS, Sh} where this list is probably far from being exhaustive).

At the turn of the millennium, Coates and many others initiated a major study towards the formulation a non-commutative main conjecture (see \cite{Bu09, Bu15, BV, CF+, CDLSS, FK, Kak, V05, V07}) which leads to a plethora of new problems and conjectures. In this context, the characteristic element now lives in a localized $K_1$-group. The goal of this paper is to introduce and make sense of this characteristic element modulo $p$. We shall describe this briefly in the introduction leaving the details to the body of the paper.

Let $G$ be a compact $p$-adic Lie group with no $p$-torsion, which is further assumed to contain a closed normal subgroup $H$ such that $G/H\cong \Zp$. Write $\pi=\pi_G$ for the natural projection map
\[ \Zp\ps{G}\twoheadrightarrow \Fp\ps{G}.\]
Following \cite{CF+, V05}, we define
\[\Sigma = \big\{x\in \Zp\ps{G}~\big|~\Zp\ps{G}/\Zp\ps{G}x \mbox{~is finitely generated over~} \Zp\ps{H} \big\},\]
where this latter set has been shown by Coates \textit{et al} to be an Ore set. As a consequence, we may speak of the localized ring $\Zp\ps{G}_{\Sigma}$. In the non-commutative Iwasawa theoretical context, the characteristic element of Iwasawa theoretical object of interest usually lives in $K_1(\Zp\ps{G}_{\Sigma})$ under appropriate assumptions. 

We now set
\[\overline{\Sigma} = \big\{z\in \Fp\ps{G}~\big|~\Fp\ps{G}/\Fp\ps{G}z \mbox{~is finitely generated over~} \Fp\ps{H} \big\}.\]
Via a similar argument of Coates \textit{et al}, one can show that $\overline{\Sigma}$ is an Ore set for $\Fp\ps{G}$. Therefore, similarly as before, one may consider the localized ring $\Fp\ps{G}_{\overline{\Sigma}}$. As will be seen in the body of the paper, it can be shown that $\pi(\Sigma) = \overline{\Sigma}$ (see Lemma \ref{Ore}). This in turn implies that the map $\pi$ induces a homomorphism of the localized rings
\[ \pi':\Zp\ps{G}_{\Sigma}\lra \Fp\ps{G}_{\overline{\Sigma}} \]
which in turn induces a group homomorphism
\[ \pi_1: K_1(\Zp\ps{G}_{\Sigma})\lra K_1(\Fp\ps{G}_{\overline{\Sigma}}).  \]
This map $\pi_1$ will play the role for our proposed notion of modulo $p$. We will make some preliminary study on this map and its properties in the body of the paper. 

As an application of the above developed theory, we study the Greenberg Selmer group attached to the tensor product of $p$-ordinary primitive Hecke eigenforms $f_1, f_2,..., f_t,h$, where $p$ is an Eisenstein prime of $h$. For simplicity, in this introduction, we shall assume that the coefficients of these eigenforms lie in $\Q$, and that the eigenforms have trivial nebentypus.  By the Eisenstein hypothesis, there is an exact sequence of $\Gal(\overline{\Q}/\Q)$-modules
\[ 0\lra \overline{C}_1 \lra A_h[\varpi]\lra \overline{C}_2 \lra  0, \]
where $\overline{C}_1$ and $\overline{C}_2$ are (free) $\Fp$-modules of dimension one. For $i=1,2$, if $\overline{\xi}_i: \Gal(\overline{\Q}/\Q) \lra \Fp^\times$ is the character attached to $\overline{C}_i$, we choose a Teichm\"uller lift of $\xi_i$ such that $\xi_i$ is a lift of $\overline{\xi}_i$ (see \cite[Lemma 2.1]{JSV}). 

Let $X_{Gr}(A_{\underline{f}\ot h}/F_\infty)$ denote the Pontryagin dual of the Greenberg Selmer group associated to the tensor product of modular forms $f_1,...,f_t, h$, and $X_{Gr}(A_{\underline{f}(\xi_i)}/F_\infty)$ the Pontryagin dual of the Greenberg Selmer group associated to the tensor product of modular forms $f_1,...,f_t$ with $C_i$ (see body of paper for their precise definitions). Under appropriate assumptions, we can attach characteristic elements to these dual Selmer groups following a procedure of Burns-Venjakob \cite{Bu09, BV} which is in turn a refinement of that in the five authors' paper \cite{CF+}. 

Our main result is the following which compares the characteristic elements of the Selmer group modulo $p$ in the sense of the above discussion (also see body of paper).

\bt
Under appropriate assumptions (see Theorem \ref{main} and Corollary \ref{maincor} for details), we have
\[ \pi_1\big(\ch_{G,\ga}\big(X_{Gr}(A_{\underline{f}\ot h}/F_\infty)\big)\big) = \pi_1\big(\ch_{G,\ga}\big(X_{Gr}(A_{\underline{f}(\xi_1)}/F_\infty)\big)\big) \pi_1\big(\ch_{G,\ga}\big(X_{Gr}(A_{\underline{f}(\xi_2)}/F_\infty)\big)\big)  \]
\et

Over the cyclotomic $\Zp$-extension of $\Q$, this result was established by Jha-Shekhar-Vangala \cite{JSV}. Our result is therefore a non-commutative generalization of their result.

We should mention that we only consider the algebraic aspects for our modulo $p$ theory in this paper. But it seems very likely that our work should have some relation with the works of Coates et al \cite{CDLSS}, and the recent works of Delbourgo and his student Gilmore \cite{DelPLMS, Del, DelGil} which considered the analytic aspects. We hope to come back in subsequent works to study this in more depth.

We now conclude this introductory section with an outline of the paper. In Section \ref{alg prelim}, we begin with certain algebraic preliminaries on the Iwasawa algebra. Once we set up these basic notion, we will develop the modulo $p$ notion for element in the localized $K_1$-group  $K_1(\Zp\ps{G}_{\Sigma})$ and study some of its properties. In Section \ref{Selmer section}, we introduce the Selmer groups of Greenberg and collect some of their properties. Finally, Section \ref{Tensor section} is where we apply our modulo $p$-theory to study the Selmer group attached to a tensor product of  $p$-ordinary eigenforms with one of them being Eisenstein at $p$.

\subsection*{Acknowledgments}

The authors are very grateful to Daniel Delbourgo for his interest and comments on the subject of the paper. The authors also liked to acknowledge the late John Coates for the numerous mathematical discussions and advices (either in person or via e-mails) over the years on noncommutative Iwasawa theory which were immensely insightful and inspirational. He shall be always remembered and greatly missed.

\section{Localized K-group} \label{alg prelim}

Throughout the paper, $p$ will always denote an odd prime. Let $\Op$ be the ring of integers of some finite extension of $\Qp$. Fix a local parameter $\varpi$ for $\Op$ and write $k$ for the residue field $\Op/\varpi$.
For a given compact $p$-adic Lie group $G$, its completed group algebra over $\Op$ is then defined by
 \[ \Op\ps{G} := \plim_U \Op[G/U], \]
where $U$ runs over the open normal subgroups of $G$ and the inverse
limit is taken with respect to the canonical projection maps. A theorem of Venjakob \cite[Theorem 3.26]{V02} asserts that $\Op\ps{G}$ is a Noetherian Auslander regular ring. Furthermore, in the event that $G$ is pro-$p$ and has no $p$-torsion, the completed group algebra $\Op\ps{G}$ is known to be a local ring with no zero divisors (see \cite{Neu}).
Consequently, it admits a skew field $Q(G)$ which is flat
over $\Op\ps{G}$ (see \cite[Chapters 6 and 10]{GW} or \cite[Chapter
4, \S 9 and \S 10]{Lam}). In view of this, one can define the $\Op\ps{G}$-rank of a finitely generated $\Zp\ps{G}$-module $M$ via
\[\rank_{\Op\ps{G}}(M)  := \dim_{Q(G)} (Q(G)\ot_{\Op\ps{G}}M).\]
We shall say that $M$ is a torsion $\Op\ps{G}$-module when $\rank_{\Op\ps{G}} (M) = 0$.
It is a standard fact that $M$ is a torsion $\Op\ps{G}$-module if and only if $\Hom_{\Op\ps{G}}(M,\Op\ps{G})=0$ (for instance, see \cite[Lemma 4.2]{LimFine}). If the torsion $\Op\ps{G}$-module $M$ also satisfies $\Ext^1_{\Op\ps{G}}(M,\Op\ps{G})=0$, we say that $M$ is a pseudo-null $\Op\ps{G}$-module.

We now extend the notion of torsion modules and pseudo-null modules to the case when $G$ is a compact $p$-adic Lie group which is not necessarily pro-$p$ or has  $p$-torsion. By a well-known theorem of Lazard, the $p$-adic Lie group $G$ contains an open normal subgroup $G_0$ which is pro-$p$ with no $p$-torsion. On the other hand, it follows from \cite[Proposition 5.4.17]{NSW} that
\[\Ext^i_{\Op\ps{G}}(M,\Op\ps{G}) \cong \Ext^i_{\Op\ps{G_0}}(M,\Op\ps{G_0})\]
for every finitely generated $\Op\ps{G}$-module $M$. In view of these observations, we shall say that $M$ is a torsion $\Op\ps{G}$-module (resp., pseudo-null $\Zp\ps{G}$-module) if $\Hom_{\Op\ps{G}}(M,\Op\ps{G})=0$ (resp., $\Ext^i_{\Op\ps{G}}(M,\Zp\ps{G})=0$ for $i=0,1$). Equivalently, this is saying that $M$ is a torsion $\Op\ps{G}$-module (resp., pseudo-null $\Op\ps{G}$-module), if $M$ is a torsion $\Op\ps{G_0}$-module (resp., pseudo-null $\Op\ps{G_0}$-module) as in the preceding paragraph. (Also, compare with \cite[Discussion after Definition 2.6]{V02}).

\subsection{Characteristic element in Iwasawa theory}

Following \cite{Bu09, Bu15, BV, CF+, FK, V05, V07}, we describe the procedure of attaching characteristic elements to a certain class of $\Op\ps{G}$-modules. From now on, $G$ will always denote a compact $p$-adic Lie group which contains a closed normal subgroup $H$ such that $\Ga:= G/H \cong \Zp$.  Write $\N_H(G)$ for the category of finitely generated $\Op\ps{G}$-module $M$ satisfying the property that $M$ is finitely generated over $\Op\ps{H}$. 
To facilitate subsequent discussion, we recall an equivalent description of this class of modules as given in \cite{CF+}. Set
\[\Si : = \Si_{G} := \big\{ x\in \Op\ps{G}~ \big|~  \Op\ps{G}/\Op\ps{G}x \mbox{ is finitely generated over } \Op\ps{H} \big\}.\]
Thanks to \cite[Theorem 2.4]{CF+}, we now know that $\Si$ is a left and right Ore set of $\Op\ps{G}$. It therefore makes sense to speak of the localization of $\Op\ps{G}$ with respect to $\Si$, which in turn is denoted by $\Op\ps{G}_{\Si}$. Furthermore, it follows from \cite[Proposition 2.3]{CF+} that a finitely generated $\Op\ps{G}$-module $M$ is annihilated by $\Si$ if and only if $M$ is finitely generated over $\Op\ps{H}$. This in turn induces an identification $K_0(\Op\ps{G}, \Op\ps{G}_{\Si}) \cong K_0(\N_H(G))$.
The localization sequence in $K$-theory (for instance, see \cite{BerK}) then  yields
the following exact sequence
\[ K_1(\Op\ps{G}) \lra K_1(\Op\ps{G}_{\Si})\stackrel{\partial_G}{\lra} K_0(\Op\ps{G}, \Op\ps{G}_{\Si}) \lra K_0(\Op\ps{G}) \lra K_0(\Op\ps{G}_{\Si})\lra 0.\]
We now state the following fundamental result of the five authors \cite[Proposition 3.4]{CF+} (also see \cite[Corollary 3.8]{Wi}).

\bl
The above connecting homomorphism $\partial_G$ is surjective.
\el

In view of the above lemma, Coates \textit{et al} \cite{CF+} define a characteristic element for $M$ in $\N_H(G)$ to be an element $\xi_M\in K_1(\Zp\ps{G}_{\Si})$ such that $\partial_G(\xi_M) = -[M]$. From their definition, we see that there can be more than one choice of characteristic element. It is then natural to ask for a ``good'' choice of such a characteristic element. Pertaining to this, Burns-Venjakob \cite{BV} (also see \cite{Bu09}) proposed such a canonical choice when the group $G$ has no $p$-torsion. (In particular, see \cite[Remarks 3.2 and 3.3]{BV} and references therein for reasons why this choice is regarded canonical.) We shall describe their construction here. 

As a start, we make the following remark. Since the group $G$ has no $p$-torsion, the Iwasawa algebra $\Op\ps{G}$ has finite global dimension (see \cite[Chap.\ V, \S 2, Exercise 5]{NSW}). By \cite[Chap.\ 7, Corollary 4.3]{GW}, this in turn implies that the localized ring $\Op\ps{G}_{\Sigma}$ has finite global dimension. Hence one can identify $K_1(\Op\ps{G}_{\Sigma})$ with the group $G_1(\Op\ps{G}_{\Sigma})$ that is generated (multiplicatively) by symbols of the form $\langle\alpha~|~Y\rangle$, where $\alpha$ is an automorphism of a finitely generated $\Op\ps{G}_{\Sigma}$-module $Y$. We shall make use of this latter identification without any further mention.

Now, let $M$ be a given module in $\N_H(G)$. Fix a topological generator $\ga$ of the group $\Ga$.
Consider the completed tensor product
\[I^G_H(M)= \Op\ps{G}\widehat{\ot}_{\Op\ps{H}}M\]
which has a natural action of $\Op\ps{G}$ via multiplication on the left. With respect to this action, we define an endomorphism $\Delta_{\ga}$ of $I^G_H(M)$ by setting
\[ \Delta_{\ga}(a\ot b) = a\widetilde{\ga}^{-1}\ot \widetilde{\ga}b\]
for each $a\in \Op\ps{G}$ and $b\in M$, where $\widetilde{\ga}\in G$ is a lift of the topological generator $\ga$. We remark that $\Delta_{\ga}$ is independent of the choice of $\widetilde{\ga}$. By \cite[Lemma 2.1]{Bu09}, the map $\delta_{\gamma,M}:=\mathrm{id}_{I^G_H(M)} - \Delta_{\ga}$ induces an automorphism of the $\Op\ps{G}_{\Si}$-module $I^G_H(M)_{\Si}$ which in turn allows us to speak of
\[ \ch_{G,\gamma}(M) := \Big\langle \delta_{\gamma,M} ~\Big|~I^G_H(M)_{\Si} \Big\rangle\]
which lies in $K_1(\Op\ps{G}_{\Si})$. We emphasis that the definition of $\ch_{G,\gamma}(M)$ is dependent on the choice of $\ga$. However, since we always work with a fixed choice of $\ga$, we shall therefore drop $``\ga"$ from our notation.

By \cite[Proposition 3.1]{Bu09}, the assignment $[M]\mapsto \ch_{G,\ga}(M)$ satisfies multiplicative property for a short exact sequence of modules in $\N_H(G)$, and so induces a group homomorphism $K_0(\N_H(G)) \lra K_1(\Op\ps{G}_{\Si})$. Precomposing this group homomorphism with the map $K_0(\N_H(G)) \lra K_0(\N_H(G))$, $[M]\mapsto -[M]$, we obtain a map $\chi: K_0(\N_H(G)) \lra K_1(\Op\ps{G}_{\Si})$ defined by $-[M]\mapsto \ch_{G}(M)$. 

\bp \label{BVsplit}
The above defined map $\chi$ is a split map of the connecting homomorphism $\partial_G$.
In other words, one has $\partial_G(\ch_{G}(M)) = -[M]$. 
\ep

\bpf
This is essentially \cite[Theorem 4.1(i)]{Bu09} (also see \cite[Proposition 4.7(i)]{BV}). 
\epf

\br

\er

\subsection{Characteristic element modulo $p$}
We now consider a modulo-$p$ variant of the discussion in the preceding subsection. 
Set
\[\overline{\Si} := \big\{ z\in k\ps{G}~ \big|~  k\ps{G}/k\ps{G}z \mbox{ is finitely generated over } k\ps{H} \big\}.\]
A similar argument to that in \cite[Theorem 2.4]{CF+} shows that $\overline{\Si}$ is a left and right Ore set of $k\ps{G}$, and so we can speak of the localized ring $k\ps{G}_{\overline{\Si}}$. Furthermore, the argument in \cite[Proposition 2.3]{CF+} carries over to show that a finitely generated $k\ps{G}$-module $N$ is annihilated by $\overline{\Sigma}$ if and only if $N$ is finitely generated over $k\ps{H}$. Consequently, there is an identification $K_0(k\ps{G}, k\ps{G}_{\overline{\Si}}) \cong K_0(\overline{\N}_H(G))$, where $\overline{\N}_H(G)$ denotes the category of $k\ps{G}$-modules that are finitely generated over $k\ps{H}$. The localization sequence in $K$-theory then yields
the following exact sequence
\[ K_1(k\ps{G}) \lra K_1(k\ps{G}_{\overline{\Si}})\stackrel{\overline{\partial}_G}{\lra} K_0(k\ps{G}, k\ps{G}_{\overline{\Si}}).\]

We can then mimic Burns-Venjakob's construction to obtain a split map of $\overline{\partial}_G$. More precisely, for a module $N$ in $\overline{\N}_H(G)$, we consider the completed tensor product
\[\overline{I}^G_H(N)= k\ps{G}\widehat{\ot}_{k\ps{H}}N\]
which has a natural action of $k\ps{G}$ via multiplication on the left. Similarly as before, we define an endomorphism $\Delta_{\ga}$ of $\overline{I}^G_H(N)$ by setting
\[ \Delta_{\ga}(x\ot y) = x\widetilde{\ga}^{-1}\ot \widetilde{\ga}y\]
for each $x\in k\ps{G}$ and $y\in N$, where $\widetilde{\ga}\in G$ is a lift of the topological generator $\ga$. By \cite[Lemma 2.1]{Bu09} again, the map $\overline{\delta}_{\ga}:=\mathrm{id}_{\overline{I}^G_H(N)} - \Delta_{\ga}$ induces an automorphism of the $k\ps{G}_{\overline{\Si}}$-module $I^G_H(N)_{\overline{\Si}}$ which in turn allows us to define
\[ \overline{\ch}_{G,\ga}(N) := \Big\langle \overline{\delta}_{\ga} ~\Big|~\overline{I}^G_H(N)_{\Si} \Big\rangle \in K_1(k\ps{G}_{\overline{\Si}}).\]
Although the definition of $\overline{\ch}_{G,\gamma}(M)$ is dependent on the choice of $\ga$, we will drop $``\ga"$  for notational simplicity.

\bp \label{BVsplit}
The map $\overline{\chi}: K_0(\overline{\N}_H(G)) \lra K_1(k\ps{G}_{\overline{\Si}})$ defined by $-[N]\mapsto \overline{\ch}_{G}(N)$ is a split map of the connecting homomorphism $\overline{\partial}_G$.
\ep

\bpf
The argument in \cite[Theorem 4.1(i)]{Bu09} (or \cite[Proposition 4.7(i)]{BV}) carries over.
\epf

\br
In particular, it follows from the preceding proposition that the map $\overline{\partial}$ is surjective.
\er

Let $\pi: \Op\ps{G}\twoheadrightarrow k\ps{G}$ denote the ring homomorphism induced by the canonical projection map $\Op\twoheadrightarrow k$. The following lemma relates the two Ore sets.

\bl \label{Ore}
One has $\pi(\Si) = \overline{\Sigma}$.
\el

\bpf
 Let $s\in \Si$. Then by definition, the module $\Op\ps{G}/\Op\ps{G}s$ is finitely generated over $k\ps{H}$. A direct calculation yields the following observation
 \[k\ps{G}/k\ps{G}\pi(s) = \Big(\Op\ps{G}/\Op\ps{G}s\Big)\ot_{\Op\ps{G}}k\ps{G} \cong \Big(\Op\ps{G}/\Op\ps{G}s\Big)\ot_{\Op}k. \]
From which, it follows that $k\ps{G}/k\ps{G}\pi(s)$ is finitely generated over $k\ps{H}$. This proves
 $\pi(\Si) \subseteq \overline{\Sigma}$.
 
 For the other inclusion, we choose a pro-$p$ subgroup $J$ of $H$ which is normal in $G$, and write $\psi_J:\Op\ps{G}\twoheadrightarrow k\ps{G/J}$ and $\theta_J:k\ps{G}\twoheadrightarrow k\ps{G/J}$ for the natural surjections. Note that $\psi_J=\theta_J\circ \pi$. Now for each $z\in \overline{\Sigma}$, choose $s\in\Zp\ps{G}$ such that $\pi(s)=z$. Then it follows that $k\ps{G/J}/k\ps{G/J}\psi_J(s) = k\ps{G/J}/k\ps{G/J}\theta_J(z)$ is finite. By \cite[Lemma 2.1]{CF+}, this in turn implies that $s\in\Sigma$, thus establishing the reverse inclusion. 
 \epf

In view of the preceding lemma, the map $\pi: \Op\ps{G}\lra k\ps{G}$ induces a ring homomorphism
\[ \pi':\Op\ps{G}_{\Sigma}\lra k\ps{G}_{\overline{\Sigma}} \]
which fits into the following commutative diagram
\begin{equation*} \entrymodifiers={!! <0pt, .8ex>+} \SelectTips{eu}{}
\xymatrix{
     \Op\ps{G} \ar[r]^{} \ar[d]^{\pi} & \Op\ps{G}_{\Sigma}\ar[d]^{\pi'} \\
     k\ps{G} \ar[r]^{} &  k\ps{G}_{\overline{\Sigma}} 
     }
      \end{equation*}
with $\Op\ps{G}_{\Sigma}\otimes_{\Op\ps{G}}k\ps{G}\cong  k\ps{G}_{\overline{\Sigma}}$. As a consequence, we have the following short exact sequence
\begin{equation}\label{local ses}
 0\lra \Op\ps{G}_{\Sigma} \stackrel{\cdot \varpi}{\lra} \Op\ps{G}_{\Sigma}\lra k\ps{G}_{\overline{\Sigma}}\lra 0.
\end{equation}
 On the other hand, by the naturality of the localization sequence, we have the following commutative diagram
\[ \entrymodifiers={!! <0pt, .8ex>+} \SelectTips{eu}{}\xymatrix{
      K_1(\Op\ps{G})\ar[r]_{}\ar[d] &  K_1(\Op\ps{G}_{\Sigma}) \ar[r]^{\partial} \ar[d]_{\pi_1}&  K_{0}(\N_H(G)) \ar[d]_{\pi_0} \ar[r] & 0\\
       K_1(k\ps{G})\ar[r]_{} &  K_1(k\ps{G}_{\overline{\Sigma}}) \ar[r]^{\overline{\partial}} &  K_{0}(\overline{\N}_H(G)) \ar[r] & 0
     }\]
     with exact rows, where the vertical maps are induced by the map $\pi$. The maps $\pi_0$ and $\pi_1$ can be described as follows. For each module $M$ in $\N_H(G)$, one has
     \[\pi_0([M]) = [M/\varpi] - [M[\varpi]].\]
     
For the map $\pi_1$, we let $\langle \alpha ~|~ Y\rangle\in K_1(\Op\ps{G}_{\Sigma})$. The short exact sequence (\ref{local ses}) induces the following exact sequence
   \begin{equation} \label{tor low} 0\lra \Tor_1^{\Op\ps{G}_{\Sigma}}( k\ps{G}_{\overline{\Sigma}},Y) \lra Y \stackrel{\cdot \varpi}{\lra} Y \lra k\ps{G}_{\overline{\Sigma}}\otimes_{\Op\ps{G}_{\Sigma}} Y\lra 0.
\end{equation}
and $\Tor_i^{\Op\ps{G}_{\Sigma}}(k\ps{G}_{\overline{\Sigma}},Y) =0$ for $i\geq 2$. The  exact sequence (\ref{tor low}) in turn tells us that 
$\Tor_1^{\Op\ps{G}_{\Sigma}}(k\ps{G}_{\overline{\Sigma}},Y) = Y[\varpi]$ and $k\ps{G}_{\overline{\Sigma}}\otimes_{\Op\ps{G}_{\Sigma}} Y = Y/\varpi$. Thus, we have
     \[ \pi_1\big(\langle \alpha ~|~ Y\rangle\big) = \langle \alpha_0 ~|~ Y/\varpi\rangle\cdot\langle \alpha_1 ~|~ Y[\varpi]\rangle^{-1}, \]
     where the maps $\alpha_0$ and $\alpha_1$ are induced by $\alpha$ via the functoriality property of Tor.

We now prove the following.

\bp \label{chicommute}
The following diagram
\[ \entrymodifiers={!! <0pt, .8ex>+} \SelectTips{eu}{}
\xymatrix{
       K_{0}(\N_H(G)) \ar[r]^{\chi} \ar[d]_{\pi_0} & K_1(\Op\ps{G}_{\Sigma}) \ar[d]_{\pi_1} \\
         K_{0}(\overline{\N}_H(G)) \ar[r]^{\overline{\chi}} & K_1(k\ps{G}_{\overline{\Sigma}}) 
     }\]
     commutes.
\ep

\bpf
Let $M$ be a module in $\N_H(G)$. Then a direct computation shows that
\[ \overline{\chi}\circ\pi_0(-[M]) =  \Big\langle \delta_{\ga,\overline{I}^G_H(M/\varpi)_{\overline{\Si}}} ~\Big|~\overline{I}^G_H(M/\varpi)_{\overline{\Si}} \Big\rangle \cdot\Big\langle \delta_{\ga,\overline{I}^G_H(M[\varpi])_{\overline{\Si}}} ~\Big|~\overline{I}^G_H(M[\varpi])_{\overline{\Si}} \Big\rangle^{-1} .\]
On the other hand, we have
\begin{eqnarray*}
\pi_1\circ\chi(-[M]) \!\!\!\!&=&\!\!\!\! \pi_1\Big(\Big\langle \delta_{\ga,I^G_H(M)_{\Si}} ~\Big|~I^G_H(M)_{\Si} \Big\rangle \Big)\\
 &=&\!\!\!\! \Big\langle \delta_{\ga,I^G_H(M)_{\Si}}' ~\Big|~I^G_H(M)_{\Si}/\varpi \Big\rangle \cdot\Big\langle \delta_{\ga,I^G_H(M)_{\Si}}'' ~\Big|~I^G_H(M)_{\Si}[\varpi] \Big\rangle^{-1} 
\end{eqnarray*}
It suffices to prove that $\overline{I}^G_H(M/\varpi)_{\overline{\Si}}\cong I^G_H(M)_{\Si}/\varpi$ and $\overline{I}^G_H(M[\varpi])_{\overline{\Si}}\cong I^G_H(M)_{\Si}[\varpi]$. 
First, note that since $\Op\ps{G}$ is projective as a compact $\Op\ps{H}$-module and $\Op\ps{G}_{\Sigma}$ is flat as a $\Op\ps{G}$-module, the functor $I^G_H(-)_{\Si}$ is exact. Therefore, upon applying $I^G_H(-)_{\Si}$ to the exact sequence
\[ 0\lra M[\varpi]\lra M\stackrel{\cdot \varpi}{\lra} M \lra M/\varpi\lra 0, \]
we obtain
\[ I^G_H(M)_{\Sigma}[\varpi]\cong I^G_H(M[\varpi])_{\Sigma},\]
\[ I^G_H(M)_{\Sigma}/\varpi\cong I^G_H(M/\varpi)_{\Sigma},\]
where $M[\varpi]$ and $M/\varpi$ both lie in $\overline{\N}_H(G)$. It therefore suffices to show that for every module $N$ in $\overline{\N}_H(G)$, one has
\[ I^G_H(N)_{\Sigma} \cong \overline{I}^G_H(N)_{\overline{\Sigma}}.\]
Indeed, applying $-\otimes_{\Op\ps{H}}N$ to the exact sequence (\ref{local ses}), we obtain
\[ \Op\ps{G}_{\Sigma}\otimes_{\Op\ps{H}}N \stackrel{\cdot \varpi}{\lra} \Op\ps{G}_{\Sigma}\otimes_{\Op\ps{H}}N \lra k\ps{G}_{\overline{\Sigma}}\otimes_{\Op\ps{H}}N \lra 0. \] 
Since $N$ is annihilated by $\varpi$, the above yields 
\begin{equation} \label{local mod p} 
\Op\ps{G}_{\Sigma}\otimes_{\Op\ps{H}}N \cong k\ps{G}_{\overline{\Sigma}}\otimes_{\Op\ps{H}}N.
\end{equation}
Now since $k\ps{G}_{\overline{\Sigma}}$ and $N$ are annihilated by $\varpi$, the latter tensor product may be taken over $k\ps{H}$. This completes the proof of the proposition. 
\epf

\bc
Let $N$ be a module in $\overline{\N}_H(G)$. Viewing $N$ as a module in $\N_H(G)$ via the projection map $\Op\ps{G}\twoheadrightarrow k\ps{G}$, we have
\[\pi_1\big(\ch_{G}(N)\big) = 0.\]
\ec

\bpf
Indeed it follows from the preceding proposition that 
\[\pi_1\big(\ch_{G}(N)\big) = \pi_1\circ\chi([N]) = \overline{\chi}\circ \pi_0([N]) = \overline{\chi}([N]-[N])= 0.\]
\epf

As will be seen in the next section, the module $M$ that we are interested in usually has the property that $M[\varpi]=0$ (see Theorem \ref{Selmer[p] = 0} below). In the event of such, Proposition \ref{chicommute} yields the following.

\bc \label{chicommutecoro}
Suppose that $M$ is a module in $\N_H(G)$ with $M[\varpi]=0$. Then one has
\[\pi_1\big(\ch_{G}(M)\big) = \overline{\ch}_{G}(M/\varpi).\]
\ec

\subsection{Rank}

In this section, we describe how the characteristic element is related to $\Op\ps{H'}$-rank for a pro-$p$ open normal subgroup $H'$ of $H$. It has been long
observed in literature that this $\Op\ps{H'}$-rank serves as a higher analog of the classical $\la$-invariant (see \cite{How}). 
We continue to retain notation from previous subsection. For every pro-$p$ open normal subgroup $H'$ of $H$, 
we define two maps
\[ r_{H'} : K_0(\Op\ps{H}) \lra \Z,  \quad [M]\mapsto \rank_{\Op\ps{H'}}(M), \]
\[ s_{H'} : K_0(k\ps{H}) \lra \Z,  \quad [N]\mapsto \rank_{k\ps{H'}}(N). \]
These two maps fit into the following diagram
\[ \entrymodifiers={!! <0pt, .8ex>+} \SelectTips{eu}{}\xymatrix{
    K_1(\Op\ps{G}_{\Sigma}) \ar[r]^{\partial} \ar[d]_{\pi_1}  & K_{0}(\N_H(G)) \ar[r]^{For}\ar[d]_{\pi_0} &   K_0(\Op\ps{H}) \ar[r]^(.65){r_{H'}} \ar[d]& \Z \ar@{=}[d] \\
   K_1(k\ps{G}_{\overline{\Sigma}}) \ar[r]^{\overline{\partial}} &  K_{0}(\overline{\N}_H(G))\ar[r]^{For} &   K_0(k\ps{H})\ar[r]^(.65){s_{H'}} & \Z
     }\]
where the vertical maps are induced by the natural projection map $\Op\twoheadrightarrow k$, the two horizontal maps $For$ are the forgetful maps and the rightmost square is commutative by the fact that 
\[ \rank_{\Op\ps{H'}}(M)=   \rank_{k\ps{H'}}(M/\varpi) -  \rank_{k\ps{H'}}(M[\varpi])\]
(see \cite[Corollary 1.10]{How} or \cite[Proposition 4.12]{LimFine}).
By abuse of notation, we still write $r_{H'}$ (resp., $s_{H'}$) for the map $K_1(\Op\ps{G}_{\Sigma})\lra  \Z$ (resp., $K_1(k\ps{G}_{\overline{\Sigma}})\lra  \Z$) given by the composition of the top (resp., bottom) horizontal maps.
Combining the above discussion with Corollary \ref{chicommutecoro}, we obtain the following.

\bp 
Suppose that $M$ is a module in $\N_H(G)$ with $M[\varpi]=0$. Then one has
\[r_{H'}\big(\pi_1\big(\ch_{G}(M)\big)\big) = s_{H'}\big(\overline{\ch}_{G}(M/\varpi)\big).\]
\ep

\section{Selmer group} \label{Selmer section}

In this section, we define the Selmer groups associated to certain datum in the sense of Greenberg \cite{G89} and study some of their properties. 
Let $F$ be a number field. As before, $\Op$ will denote the ring of integers of some finite extension of $\Qp$ with local parameter $\varpi$ and residue field $k$ ($=\Op/\varpi)$.

We now introduce the axiomatic conditions on our datum which is denoted by $\big(B, \{B_v\}_{v|p}\big)$ and satisfies all of the following four conditions \textbf{(C1)-(C4)}.

\begin{enumerate}
 \item[(\textbf{C1})] $B$ is a
cofinitely generated cofree $\Op$-module of $\Op$-corank $d$ with a
continuous, $\Op$-linear $\Gal(\bar{F}/F)$-action which is
unramified outside a finite set of primes of $F$.

 \item[(\textbf{C2})] For each prime $v$ of $F$ above $p$, there is a
distinguished $\Gal(\bar{F}_v/F_v)$-submodule $B_v$ of $B$ which is cofree of
$\Op$-corank $d_v$. We shall write $B_v^- = B/B_v$ which is a cofree $\Op$-module of
$\Op$-corank $d-d_v$.

 \item[(\textbf{C3})] For each real prime $v$ of $F$, $B_v^+:=
B^{\Gal(\bar{F}_v/F_v)}$  is cofree of
$\Op$-corank $d^+_v$.

\item[(\textbf{C4})] The following equality
  \[\sum_{v|p} (d-d_v)[F_v:\Qp] = dr_2(F) +
 \sum_{v~\mathrm{real}}(d-d^+_v)\]
holds. Here $r_2(F)$ denotes the number of complex primes of $F$.
\end{enumerate}

As we need to work with Selmer groups defined over a tower of number fields, we need to consider the base change of our datum which we now do. For a finite extension $L$ of $F$, the base change of the datum $\big(B,
\{B_w\}_{w|p} \big)$ over $L$ is given as follows:

(1) $B$ can be viewed as a $\Gal(\bar{F}/L)$-module via restriction of the Galois action.

(2) For each prime $w$
of $L$ above $p$, we set $B_w =B_v$, where $v$ is the prime of $F$
below $w$, and view it as a $\Gal(\bar{F}_v/L_w)$-module via the appropriate restriction.
Note that we then have $d_w= d_v$.

(3) For each real prime $w$ of $L$ which lies above a real prime $v$ of $F$, we set
$B_w^+= B^{\Gal(\bar{F}_v/F_v)}$ and write $d^+_w
= d^+_v$.

The following lemma gives some
sufficient conditions for equality in \textbf{(C4)} to hold for
the datum $\big(B, \{B_w\}_{w|p}\big)$ over $L$.

\bl \label{data base change} Let $\big(B, \{B_v\}_{v|p}\big)$ be a datum defined over $F$.
 Suppose that at least one of the following statements holds.
 \begin{enumerate}
\item[$(i)$] $[L:F]$ is odd.

\item[$(ii)$] $F$ has no real primes.

\item[$(iii)$] $F$ is totally real, $L$ is totally imaginary and 
\[\sum_{v~\mathrm{real}}d_v^+ = d[F:\Q]/2.\]
 \end{enumerate}
Then the datum $\big(B, \{B_w\}_{w|p}\big)$ obtained via base changing satisfies \textbf{(C1)-(C4)}. In particular, we have the equality
 \[ \sum_{w|p} (d-d_w)[L_w:\Qp] = dr_2(L) +
 \sum_{w~\mathrm{real}}(d-d^+_w).\]
 \el

\bpf See \cite[Lemma 3.0.1]{LimCMu}. \epf

Let $S$ be a finite set of primes of $F$ which contains all the primes above $p$, the ramified primes of $B$ and all the infinite primes of $F$. Denote by $F_S$ the maximal algebraic extension of $F$ which is unramified outside $S$. Write $G_S({\mathcal{F}}) = \Gal(F_S/\mathcal{F})$ for every algebraic extension $\mathcal{F}$ of $F$ contained in $F_S$. 
We say that $F_{\infty}$ is an
($S$-)\textit{admissible $p$-adic Lie extension} of $F$ if (i)
$\Gal(F_{\infty}/F)$ is compact $p$-adic Lie group, (ii)
$F_{\infty}$ contains the cyclotomic $\Zp$ extension $F_{\cyc}$ of
$F$ and (iii) $F_{\infty}$ is unramified outside $S$. Write $G =
\Gal(F_{\infty}/F)$, $H = \Gal(F_{\infty}/F_{\cyc})$ and $\Ga
=\Gal(F_{\cyc}/F)$. 

For the remainder of the section, we will always assume that for every finite extension $L$ of $F$ contained in $F_\infty$, the datum $\big(B, \{B_w\}_{w|p}\big)$ obtained by base change satisfies \textbf{(C1)-(C4)}.

Let $L$ be a finite extension of $F$ contained in $F_S$. Consider the following local condition
\[ H^1_{Gr}(L_w, B)=
\begin{cases} \ker\big(H^1(L_w, B)\lra H^1(L_w^{ur}, B^-_w)\big) & \text{\mbox{if} $w|p$},\\
 \ker\big(H^1(L_w, B)\lra H^1(L^{ur}_w, B)\big) & \text{\mbox{if} $w\nmid p$,}
\end{cases} \]
where $L_w^{ur}$ is the maximal unramified extension of $L_w$.
The Greenberg Selmer group attached to the datum $\big(B,
\{B_w\}_{w|p} \big)$  is then defined by
\[ \Sel_{Gr}(B/L) = \ker\left( H^1(G_S(L),B)\lra \bigoplus_{v \in S}\bigoplus_{w|v}\frac{H^1(L_w, B)}{H^1_{Gr}(L_w,
B)}\right).\]  
We define $\Sel_{Gr}(B/F_{\infty}) = \ilim_L \Sel_{Gr}(B/L)$,
where the limit runs over all finite extensions $L$ of $F$ contained
in $F_{\infty}$. The Pontryagin dual of $\Sel_{Gr}(B/F_{\infty})$ is then denoted by $X_{Gr}(B/F_{\infty})$.  

We shall introduce another Selmer group of Greenberg which is called the strict Selmer group. This group is useful to work with at times which we will do in the paper. Set
\[ H^1_{str}(L_w, B)=
\begin{cases} \ker\big(H^1(L_w, B)\lra H^1(L_w, B^-_w)\big) & \text{\mbox{if} $w|p$},\\
 \ker\big(H^1(L_w, B)\lra H^1(L^{ur}_w, B)\big) & \text{\mbox{if} $w\nmid p$.}
\end{cases} \]
Note that the local condition $H^1_{str}(L_w, B)$ is the same as $H^1_{Gr}(L_w, B)$ for primes not dividing $p$. The strict Selmer group attached to the datum $\big(B,
\{B_w\}_{w|p}\big)$ is defined to be
\[ \Sel_{str}(B/L) = \ker\left( H^1(G_S(L),B)\lra \bigoplus_{v \in S}\bigoplus_{w|v}\frac{H^1(L_w, B)}{H^1_{str}(L_w,
B)}\right).\]  

We set $\Sel_{str}(B/F_{\infty}) = \ilim_L \Sel_{str}(B/L)$,
where the limit runs over all finite extensions $L$ of $F$ contained
in $F_{\infty}$. The Pontryagin dual of $\Sel_{str}(B/F_{\infty})$ is then denoted by $X_{str}(B/F_{\infty})$.

\bl \label{localcondition}
For each $v\in S$, we fix a prime of $F_\infty$ above $v$ which we, by abuse of notation, denote by $v$. We then write $G_v$ for the decomposition group of $G$ at $v$. Then one has the following identification
\[ \ilim_L\bigoplus_{w|v}\frac{H^1(L_w, B)}{H^1_{Gr}(L_w,
B)} \cong\left\{
           \begin{array}{ll}
             \mathrm{Coind}^G_{G_v} \Big(H^1(F_{\infty, v}^{ur}, B^-_v)^{\Gal(F_{\infty,v}^{ur}/F_{\infty,v})}\Big), & \hbox{if $v\mid p$} \\
             \mathrm{Coind}^G_{G_v} \Big(H^1(F_{\infty,v},B)\Big), & \hbox{if $v\nmid p$.}
           \end{array}
         \right.
  \]

\[ \ilim_L\bigoplus_{w|v}\frac{H^1(L_w, B)}{H^1_{str}(L_w,
B)} \cong\left\{
           \begin{array}{ll}
             \mathrm{Coind}^G_{G_v} \Big(H^1(F_{\infty, v}, B^-_v)\Big), & \hbox{if $v\mid p$} \\
             \mathrm{Coind}^G_{G_v} \Big(H^1(F_{\infty,v},B)\Big), & \hbox{if $v\nmid p$.}
           \end{array}
         \right.
  \]
\el

\bpf
We shall prove the case of $v\mid p$ for the $Gr$ condition, the remaining cases can be proven similar. By definition, we have the exact sequence
\[  0 \lra H^1_{Gr}(L_v, B)\lra H^1(L_v, B) \lra H^1(L_v^{ur}, B^-_v)\]
which upon taking limit yields
\begin{equation}\label{Gr sesA}
  0 \lra \ilim H^1_{Gr}(L_v, B)\lra H^1(F_{\infty,v}, B) \lra H^1(F_{\infty, v}^{ur}, B^-_v)
\end{equation}
Note that the map $H^1(F_{\infty,v}, B) \lra H^1(F_{\infty, v}^{ur}, B^-_v)$ is given by the composition
\[H^1(F_{\infty,v}, B) \lra H^1(F_{\infty,v}, B^-_v) \lra H^1(F_{\infty, v}^{ur}, B^-_v)^{\Gal(F_{\infty,v}^{ur}/F_{\infty,v})}\subseteq H^1(F_{\infty, v}^{ur}, B^-_v), \]
where the first map is surjective by \cite[Theorem 7.1.8]{NSW} and the second map is surjective as a consequence of $\Gal(F_{\infty,v}^{ur}/F_{\infty,v})$ having $p$-cohomological dimension one. Thus, from the sequence (\ref{Gr sesA}), we obtain the short exact sequence 
\begin{equation}\label{Gr ses}
  0 \lra \ilim H^1_{Gr}(L_v, B)\lra H^1(F_{\infty,v}, B) \lra H^1(F_{\infty, v}^{ur}, B^-_v)^{\Gal(F_{\infty,v}^{ur}/F_{\infty,v})}\lra 0.
\end{equation}
This proves what we set to show.
\epf

To simplify notation, we shall write \[J_v^{Gr}(B/F_\infty)=\left\{
           \begin{array}{ll}
             \mathrm{Coind}^G_{G_v} \big(H^1(F_{\infty, v}^{ur}, B^-_v)^{\Gal(F_{\infty,v}^{ur}/F_{\infty,v})}\big), & \hbox{if $v\mid p$} \\
             \mathrm{Coind}^G_{G_v} \big(H^1(F_{\infty,v},B)\big), & \hbox{if $v\nmid p$,}
           \end{array}
         \right.\]
and
\[J_v^{str}(B/F_\infty)=\left\{
           \begin{array}{ll}
             \mathrm{Coind}^G_{G_v} \big(H^1(F_{\infty, v}^{ur}, B^-_v)\big), & \hbox{if $v\mid p$} \\
             \mathrm{Coind}^G_{G_v} \big(H^1(F_{\infty,v},B)\big), & \hbox{if $v\nmid p$.}
           \end{array}
         \right.\]

Note that $J_v^{str}(B/F_\infty)=J_v^{Gr}(B/F_\infty)$ for $v\nmid p$. Furthermore, if $v$ is not divisible by $p$ and the decomposition group of $G$ at $v$ has dimension two, then $F_{\infty,v}$ has no non-trivial $p$-extension (cf. \cite[Theorem 7.5.3]{NSW}) and so $J_v^{str}(B/F_\infty)=J_v^{Gr}(B/F_\infty)=0$.

For data coming from ordinary representations, it is expected that $X_{Gr}(B/F_\infty)$ is a torsion $\Op\ps{G}$-module (see \cite[Conjecture 1]{G89}). For our purposes, we require $X_{Gr}(B/F_\infty)$ to be in $\N_H(G)$, or equivalently, $X_{Gr}(B/F_\infty)$ is a finitely generated $\Op\ps{H}$-module. The next lemma provides a useful criterion for this finite generation property.

\bl \label{fg NHG}
Suppose that there exists a finite extension $L$ of $F$ contained in $F_\infty$ such that $X_{Gr}(B/L_\cyc)$ is finitely generated over $\Op$ and $\Gal(F_\infty/L)$ is pro-$p$. Then  $X_{Gr}(B/F_\infty)$ is finitely generated over $\Op\ps{H}$.
\el

\bpf
This can be proven by a similar argument to that in \cite[Theorem 2.1]{CS12}. 
\epf



We now define a mod-$\varpi$ variant of the Greenberg Selmer group attached to the datum $\big(B,
\{B_v\}_{v|p}\big)$. In this context, the said group is defined by
\[ \Sel_{Gr}(B[\varpi]/L) = \ker\left( H^1(G_S(L),B[\varpi])\lra \bigoplus_{v \in S}\bigoplus_{w|v}\frac{H^1(L_w, B[\varpi])}{H^1_{Gr}(L_w,
B[\varpi])}\right).\]  

We define $\Sel_{Gr}(B[\varpi]/F_{\infty}) = \ilim_L \Sel_{Gr}(B[\varpi]/L)$,
where the limit runs over all finite extensions $L$ of $F$ contained
in $F_{\infty}$. 
The mod-$\varpi$ strict Selmer group $\Sel_{str}(B[\varpi]/F_{\infty})$ is defined via the same modification.
We write $X_{Gr}(B[\varpi]/F_{\infty})$ and $X_{str}(B[\varpi]/F_{\infty})$ for the Pontryagin dual of $\Sel_{Gr}(B[\varpi]/F_{\infty})$ and $\Sel_{str}(B[\varpi]/F_{\infty})$ respectively.  

One can obtain a similar result to that in Lemma \ref{localcondition}, and as before, we shall write 

\[J_v^{Gr}(B[\varpi]/F_\infty)=\left\{
           \begin{array}{ll}
             \mathrm{Coind}^G_{G_v} \Big(H^1(F_{\infty,v}^{ur},B^-_v[\varpi])^{\Gal(F_{\infty,v}^{ur}/F_{\infty,v})}\Big), & \hbox{if $v\mid p$} \\
             \mathrm{Coind}^G_{G_v}\Big(H^1(F_{\infty,v},B[\varpi])\Big), & \hbox{if $v\nmid p$.}
           \end{array}
         \right.\]
and
\[J_v^{str}(B/F_\infty)=\left\{
           \begin{array}{ll}
             \mathrm{Coind}^G_{G_v} \Big(H^1(F_{\infty, v}^{ur}, B^-_v[\varpi])\Big), & \hbox{if $v\mid p$} \\
             \mathrm{Coind}^G_{G_v} \Big(H^1(F_{\infty,v},B[\varpi])\Big), & \hbox{if $v\nmid p$.}
           \end{array}
         \right.\]

We record a preliminary lemma.

\bl \label{fg H}
The following statements are equivalent. 
\begin{enumerate}
  \item[$(a)$] $X_{Gr}(B/F_\infty)$ is finitely generated over $\Op\ps{H}$.
  \item[$(b)$]  $X_{Gr}(B[\varpi]/F_{\infty})$ is finitely generated over $k\ps{H}$.
  \item[$(c)$] $X_{str}(B/F_\infty)$ is finitely generated over $\Op\ps{H}$.
  \item[$(d)$] $X_{str}(B[\varpi]/F_\infty)$ is finitely generated over $k\ps{H}$.
\end{enumerate}
\el

\bpf
Consider the following commutative diagram
\[ \entrymodifiers={!! <0pt, .8ex>+} \SelectTips{eu}{}\xymatrix{
     0 \ar[r] &  \Sel_{Gr}(B[\varpi]/F_\infty)\ar[r] \ar[d]_{s} &  H^1(G_S(F_\infty),B[\varpi]) \ar[r] \ar[d]_{h}&  \bigoplus_{v\in S} J_v^{Gr}(B[\varpi]/F_\infty) \ar[d]_{g=\oplus_v g_v}\\
     0 \ar[r] &  \Sel_{Gr}(B/F_\infty)[\varpi]\ar[r]_{} &   H^1(G_S(F_\infty),B)[\varpi] \ar[r] &  \bigoplus_{v\in S}J_v^{Gr}(B/F_\infty)[\varpi]
     }\]
with exact rows.
     The long exact sequence in cohomology arising from the short exact sequence
     \[ 0 \lra B[\varpi] \lra B \stackrel{\cdot\varpi}{\lra} B \lra 0 \]
shows that the map $h$ is surjective with  $\ker h = B(F_\infty)/\varpi$ which is cofinitely generated over $\Op$.
Similarly, for each $v\in S$, one sees that $\ker g_v$ is equal to $\mathrm{Coind}^{G}_{G_v}\big(B^-_v(F_{\infty,v})\big)/\varpi$ or $\mathrm{Coind}^{G}_{G_v}\Big(B(F_{\infty,v})/\varpi\Big)$
accordingly as $v$ does or does not divide $p$. Since $\mathrm{Coind}^{G}_{G_v}(-)= \oplus_{w\mid v}\mathrm{Coind}^{H}_{H_w}(-)$ and there are only finitely many primes of $F_\cyc$ above $v$, we see that $\ker g_v$ is cofinitely generated over $k\ps{H}$. Consequently, the map $s$ has kernel and cokernel which are cofinitely generated over $k\ps{H}$. The equivalence of (a) and (b) is now immediate from this. Applying a similar argument to the strict Selmer group, we also have the equivalence between (c) and (d).

On the other hand, the defining sequences of the mod-$\varpi$ Selmer groups fit into the following commutative diagram
\[ \entrymodifiers={!! <0pt, .8ex>+} \SelectTips{eu}{}\xymatrix{
     0 \ar[r] &  \Sel_{str}(B[\varpi]/F_\infty)\ar[r] \ar[d]_{} &  H^1(G_S(F_\infty),B[\varpi]) \ar[r] \ar@{=}[d] &  \bigoplus_{v\in S} J_v^{str}(B[\varpi]/F_\infty) \ar[d]_{\alpha=\oplus_v \alpha_v}\\
     0 \ar[r] &  \Sel_{Gr}(B[\varpi]/F_\infty)\ar[r]_{} &   H^1(G_S(F_\infty),B[\varpi]) \ar[r] &  \bigoplus_{v\in S}J_v^{Gr}(B[\varpi]/F_\infty)
     }\]
with exact rows. Plainly $\alpha_v$ is the identity map whenever $v\nmid p$. When $v\mid p$, the inflation-restriction sequence yields the following exact sequence
\begin{multline*}
  0\lra H^1(\Gal(F_{\infty,v}^{ur}/F_{\infty,v}),  B^-_v(F_{\infty,v}^{ur})[\varpi])) \lra H^1(F_{\infty,v}, B^-_v[\varpi]) \\
 \lra H^1(F_{\infty,v}^{ur},B_v^-[\varpi])^{\Gal(F_{\infty,v}^{ur}/F_{\infty,v})} \lra 0
\end{multline*}
which in turn implies that 
\[\ker\alpha_v = \mathrm{Coind}^{G}_{G_v}H^1(\Gal(F_{\infty,v}^{ur}/F_{\infty,v}),  B^-_v(F_{\infty,v}^{ur})[\varpi])). \]  
Now taking these observations into consideration, upon applying a snake lemma argument, we obtain
\begin{multline*}
  \hspace{.5in} 0 \lra \Sel_{str}(B[\varpi]/F_\infty) \lra \Sel_{Gr}(B[\varpi]/F_\infty) \\
  \lra \bigoplus_{v\mid p}\mathrm{Coind}^{G}_{G_v}H^1(\Gal(F_{\infty,v}^{ur}/F_{\infty,v}),  B^-_v(F_{\infty,v}^{ur})[\varpi])) \hspace{.5in}
\end{multline*}
As there are only finitely many primes of $F_\cyc$ above $v$, we see that the rightmost term is cofinitely generated over $k\ps{H}$. This therefore establishes the equivalence of (b) and (d). 
\epf

\bp \label{p-adic torsion}
Let $F_\infty$ be an admissible $p$-adic Lie extension of $F$. Suppose that $X_{str}(B/F_\infty)$ is finitely generated over $\Op\ps{H}$. Then the following assertions are valid.

\begin{enumerate}
\item[$(a)$] $H^2(G_S(F_\infty),B[\varpi])=0$.
\item[$(b)$] There is a short exact sequence
\[ 0\lra  \Sel_{str}(B[\varpi]/F_\infty)\lra H^1(G_S(F_\infty),B[\varpi])\lra \bigoplus_{v\in S}J_v^{str}(B[\varpi]/F_\infty)\lra 0.\]
\item[$(c)$] There is a short exact sequence
\[ 0\lra  \Sel_{Gr}(B[\varpi]/F_\infty)\lra H^1(G_S(F_\infty),B[\varpi])\lra \bigoplus_{v\in S}J_v^{Gr}(B[\varpi]/F_\infty)\lra 0.\]
\item[$(d)$] One has the following short exact sequence \begin{multline*}
  \hspace{.5in} 0 \lra \Sel_{str}(B[\varpi]/F_\infty) \lra \Sel_{Gr}(B[\varpi]/F_\infty) \\
  \lra \bigoplus_{v\mid p}\mathrm{Coind}^{G}_{G_v}H^1(\Gal(F_{\infty,v}^{ur}/F_{\infty,v}),  B_v^-(F_{\infty,v}^{ur})[\varpi])) \lra 0\hspace{.5in}
\end{multline*}
\end{enumerate}
\ep

\bpf
Combining the hypothesis of the proposition with Lemma \ref{fg H}, we see that $\Sel_{Gr}(B[\varpi]/F_\infty)$ is cofinitely generated over $k\ps{H}$. In particular, 
$\Sel_{Gr}(B[\varpi]/F_\infty)$ is a cotorsion $k\ps{G}$-module.

Now set $B^*= \Hom_{\Zp}(\plim_n B[\varpi^n], \mu_{p^\infty})$ and define $H^1_{\mathrm{Iw}}(F_\infty/F, B^*[\varpi]):=\plim_L H^1(G_S(L), B^*[\varpi])$. From the Poitou-Tate sequence, we have
\[ 0\lra \Sel_{str}(B[\varpi]/F_\infty) \lra H^1(G_S(F_\infty),B[\varpi])\lra \bigoplus_{v\in S} J_v^{str}(B[\varpi]/F_\infty)\]
\[ \lra \mathfrak{S}_{str}(B[\varpi]/F_\infty)^\vee \lra H^2(G_S(F_\infty),B[\varpi]) \lra 0, \]
where $\mathfrak{S}_{str}(B[\varpi]/F_\infty)$ is a submodule of $H^1_{\mathrm{Iw}}(F_\infty/F, B^*[\varpi])$.
Standard corank calculations (\cite[Proposition 1-3]{G89} or \cite[Theorems 3.2 and 4.1]{OV}) tell us that
\[ \corank_{k\ps{G}}H^1(G_S(F_\infty),B[\varpi]) - \corank_{k\ps{G}}H^2(G_S(F_\infty),B[\varpi]) = dr_2(F) +
 \sum_{v~\mathrm{real}}(d-d^+_v),\]
\[ \corank_{k\ps{G}}J_v(B[\varpi]/F_\infty)=\left\{
                                             \begin{array}{ll}
                                               d_v[F_v:\Qp], & \hbox{if $v\mid p$;} \\
                                               0, & \hbox{if $v\nmid p$.}
                                             \end{array}
                                           \right.
 \]
Putting these information into the Poitou-Tate sequence and taking our standing hypothesis  $\mathbf{(C4)}$ into consideration, we have that $\mathfrak{S}_{str}(B/F_\infty)$ is torsion over $k\ps{G}$.

On the other hand, the low degree term of the spectral sequence of Jannsen \cite[Theorems 1 and 11]{jannsen}
\[ \Ext^i_{k\ps{G}}\big(H^j(G_S(F_\infty),  B^*[\varpi])^\vee,k\ps{G}\big)\Longrightarrow H^{i+j}_{\mathrm{Iw}}(F_\infty/F, B^*[\varpi])\] 
yields an exact sequence
\[ 0\lra \Ext^1_{k\ps{G}}\big(H^0(G_S(F_\infty),  B^*[\varpi])^\vee,k\ps{G}\big)\longrightarrow H^{1}_{\mathrm{Iw}}(F_\infty/F, B^*[\varpi])\]\[ \lra \Hom_{k\ps{G}}\big(H^1(G_S(F_\infty),  B^*[\varpi])^\vee,k\ps{G}\big)\]
Since $H^0(G_S(F_\infty),  B^*[\varpi])^\vee$ is finite and $\dim G\geq 2$, $H^0(G_S(F_\infty),  B^*[\varpi])^\vee$ is torsion over $k\ps{H}$ and hence by \cite{V03}, it is pseudo-null over $k\ps{G}$. Therefore, the first term in the above exact sequence vanishes which in turn implies that $H^{1}_{\mathrm{Iw}}(F_\infty/F, B^*[\varpi])$, and  $\mathfrak{S}_{str}(B[\varpi]/F_\infty)$, injecting into a $\Hom_{k\ps{G}}$-term. Thus, it follows that $\mathfrak{S}_{str}(B[\varpi]/F_\infty)$ is torsionfree over $k\ps{G}$. Combining this with the observation in the preceding paragraph, we see that $\mathfrak{S}_{str}(B[\varpi]/F_\infty)$ must be trivial. This proves statements (a) and (b). 

Now consider the following commutative diagram
\[ \entrymodifiers={!! <0pt, .8ex>+} \SelectTips{eu}{}\xymatrix{
     0 \ar[r] &  \Sel_{str}(B[\varpi]/F_\infty)\ar[r] \ar[d]_{} &  H^1(G_S(F_\infty),B[\varpi]) \ar[r] \ar@{=}[d] &  \bigoplus_{v\in S} J_v^{str}(B[\varpi]/F_\infty) \ar[d]_{\al=\oplus_v \al_v} \ar[r] & 0\\
     0 \ar[r] &  \Sel_{Gr}(B[\varpi]/F_\infty)\ar[r]_{} &   H^1(G_S(F_\infty),B[\varpi]) \ar[r] &  \bigoplus_{v\in S}J_v^{Gr}(B[\varpi]/F_\infty) &
     }\]
with exact rows, where the final zero in the top row is a consequence of (b). The rightmost vertical map can be easily checked to be surjective
by definition. Hence it follows that the rightmost map of the bottom row is surjective which yields (c). Finally, the exact sequence (d) follows from these observations and an application of the snake lemma.
\epf

By a similar argument to that in the preceding proposition, we have the following.

\bp \label{p-adic torsion2}
Let $F_\infty$ be an admissible $p$-adic Lie extension of $F$. Suppose that $X(B/F_\infty)$ is finitely generated over $\Op\ps{H}$. Then the following assertions are valid.

\begin{enumerate}
\item[$(a)$] $H^2(G_S(F_\infty),B)=0$.
\item[$(b)$] There is a short exact sequence
\[ 0\lra  \Sel_{str}(B/F_\infty)\lra H^1(G_S(F_\infty),B)\lra \bigoplus_{v\in S}J_v^{str}(B/F_\infty)\lra 0.\]
\item[$(c)$] There is a short exact sequence
\[ 0\lra  \Sel_{Gr}(B/F_\infty)\lra H^1(G_S(F_\infty),B)\lra \bigoplus_{v\in S}J_v^{Gr}(B/F_\infty)\lra 0.\]
\item[$(d)$] One has the following short exact sequence \begin{multline*}
  \hspace{.5in} 0 \lra \Sel_{str}(B/F_\infty) \lra \Sel_{Gr}(B/F_\infty) \\
  \lra \bigoplus_{v\mid p}\mathrm{Coind}^{G}_{G_v}H^1(\Gal(F_{\infty,v}^{ur}/F_{\infty,v}),  B/B_v(F_{\infty,v}^{ur}))) \lra 0\hspace{.5in}
\end{multline*}
\end{enumerate}
\ep

We now prove the following result on the structure on the Selmer groups which will allow us to apply Corollary \ref{chicommutecoro}.

\bt \label{Selmer[p] = 0}
Let $F_\infty$ be an admissible $p$-adic Lie extension of $F$. Suppose that $X_{Gr}(B/F_\infty)$ is finitely generated over $\Op\ps{H}$. Then one has $X_{Gr}(B/F_\infty)[\varpi]=0$ and  $X_{str}(B/F_\infty)[\varpi]=0$.
\et

\bpf
We will give the proof of the first equality, the second being similar.
Consider the following diagram
\begin{equation}
\entrymodifiers={!! <0pt, .8ex>+} \SelectTips{eu}{}\xymatrix{
     0 \ar[r] &  \Sel_{Gr}(B[\varpi]/F_\infty)\ar[r] \ar[d]_{s} &  H^1(G_S(F_\infty),B[\varpi]) \ar[r] \ar[d]_{h}&  \bigoplus_{v\in S} J_v^{Gr}(B[\varpi]/F_\infty) \ar[d]_{g=\oplus_v g_v} \ar[r] & 0 \\
     0 \ar[r] &  \Sel_{Gr}(B/F_\infty)[\varpi]\ar[r]_{} &   H^1(G_S(F_\infty),B)[\varpi] \ar[r] &  \bigoplus_{v\in S}J_v^{Gr}(B/F_\infty)[\varpi] \ar[r] & \cdots
     }
\end{equation}
where the bottom row is part of the following exact sequence
\begin{multline*}
  0 \lra \Sel_{Gr}(B/F_\infty)[\varpi] \lra   H^1(G_S(F_\infty),B)[\varpi] \lra  \bigoplus_{v\in S}J_v^{Gr}(B/F_\infty)[\varpi] \\
  \lra \Sel_{Gr}(B/F_\infty)/\varpi \lra H^1(G_S(F_\infty),B)/\varpi
\end{multline*}
in view of Proposition \ref{p-adic torsion2}(b). 
As the rightmost vertical map of the diagram is surjective, the bottom row is a short exact sequence and so we have an injection
\[ \Sel_{Gr}(B/F_\infty)/\varpi \hookrightarrow H^1(G_S(F_\infty),B)/\varpi.\]
But the latter injects into $H^2(G_S(F_\infty),B[\varpi])$ which vanishes by Proposition \ref{p-adic torsion}(a).
Hence it follows that $\Sel_{Gr}(B/F_\infty)/\varpi  =0$, or equivalently, $X_{Gr}(B/F_\infty)[\varpi]=0$.
\epf

We end the section with the following observation on $\ch_{G}\big(X_{str}(B/F_\infty)\big)$ and $\ch_{G}\big(X_{Gr}(B/F_\infty)\big)$. 

\bp \label{char mod p}
Let $F_\infty$ be an admissible $p$-adic Lie extension of $F$. Suppose that $X_{Gr}(B/F_\infty)$ is finitely generated over $\Op\ps{H}$. Then one has
\[ \pi_1\big(\ch_{G}\big(X_{str}(B/F_\infty)\big)\big) = \overline{\ch}_{G}\big(X_{str}(B[\varpi]/F_\infty)\big) ~\overline{\ch}_{G}\big(B(F_\infty)^\vee[\varpi]\big) \]
\[\hspace{1.5in} \times\prod_{v\in S_p}\overline{\ch}_{G}\left(\mathrm{Coind}^{G}_{G_v}\big(B^-_v(F_{\infty,v})/\varpi\big)^\vee
\right)^{-1}\]
\[\hspace{1.75in}\times\prod_{v\in S\setminus S_p}\overline{\ch}_{G}\left(\mathrm{Coind}^{G}_{G_v}\Big(B(F_{\infty,v})/\varpi\Big)^\vee\right)^{-1}\]
and 
\[ \ch_{G}\big(X_{Gr}(B/F_\infty)\big) = \ch_{G}\big(X_{str}(B/F_\infty)\big)\prod_{v\in S_p}\ch_{G}\left(\mathrm{Coind}^{G}_{G_v}H^1(\Gal(F_{\infty,v}^{ur}/F_{\infty,v}),  B^-_v(F_{\infty,v}^{ur}))^\vee\right)\]
\ep

\bpf
The second equality is immediate from Proposition \ref{p-adic torsion2}(d). For the first, note that
by virtue of Theorem \ref{Selmer[p] = 0}, we may apply Corollary \ref{chicommutecoro} to obtain
\[ \pi_1\big(\ch_{G}\big(X_{str}(B/F_\infty)\big)\big) = \overline{\ch}_{G}\big(X_{str}(B/F_\infty)/\varpi\big).\]
Taking this into account, we are reduced to showing the equality
\[\big[X_{str}(B/F_\infty)/\varpi\big] = \big[X_{str}(B[\varpi]/F_\infty)\big] + \big[B(F_\infty)^\vee[\varpi]\big] -\sum_{v\in S_p}\left[\mathrm{Coind}^{G}_{G_v}\big(B^-_v(F_{\infty,v})/\varpi\big)^\vee
\right]\] 
 \[  - \sum_{v\in S\setminus S_p}\left[\mathrm{Coind}^{G}_{G_v}\Big(B(F_{\infty,v})/\varpi\Big)^\vee\right]\]
in $K_0(\overline{\N}_H(G))$. This said equality in turn follows from an application of a snake lemma argument to the following commutative diagram
\[
\entrymodifiers={!! <0pt, .8ex>+} \SelectTips{eu}{}\xymatrix{
     0 \ar[r] &  \Sel_{str}(B[\varpi]/F_\infty)\ar[r] \ar[d]_{s} &  H^1(G_S(F_\infty),B[\varpi]) \ar[r] \ar[d]_{h}&  \bigoplus_{v\in S} J_v^{str}(B[\varpi]/F_\infty) \ar[d]_{g=\oplus_v g_v} \ar[r] & 0\\
     0 \ar[r] &  \Sel_{str}(B/F_\infty)[\varpi]\ar[r]_{} &   H^1(G_S(F_\infty),B)[\varpi] \ar[r] &  \bigoplus_{v\in S}J_v^{str}(B/F_\infty)[\varpi] &
     }
\]
with exact rows.
\epf

\bp \label{char mod p rank}
Retain the assumptions of Proposition \ref{char mod p}. Then for every open normal pro-$p$ subgroup $H'$ of $H$, we have
\begin{multline*}
\rank_{\Op\ps{H'}}\big(X_{str}(B/F_\infty)\big) = \rank_{k\ps{H'}}\big(X_{str}(B[\varpi]/F_\infty)\big)  -\sum_{v\in S_p}\rank_{k\ps{H'}}\left(\mathrm{Coind}^{G}_{G_v}\big(B^-_v(F_{\infty,v})/\varpi\big)^\vee
\right)\\
-\sum_{v\in S\setminus S_p}\rank_{k\ps{H'}}\left(\mathrm{Coind}^{G}_{G_v}\Big(B(F_{\infty,v})/\varpi\Big)^\vee\right)
\end{multline*}
and 
\begin{multline*}
  \rank_{\Op\ps{H'}}\big(X_{Gr}(B/F_\infty)\big) =  \rank_{\Op\ps{H'}}\big(X_{str}(B/F_\infty)\big)  \\
  +\sum_{v\in S_p}\rank_{\Op\ps{H'}}\left(\mathrm{Coind}^{G}_{G_v}H^1(\Gal(F_{\infty,v}^{ur}/F_{\infty,v}),  B^-_v(F_{\infty,v}^{ur}))^\vee\right)
\end{multline*}
\ep

\section{Tensor product of modular forms} \label{Tensor section}

From now on, $t$ will denote a fixed positive integer.
Suppose that we are given $p$-ordinary primitive Hecke eigenforms $f_1, f_2,..., f_t,h$ with weights $k_1, k_2,...,k_t, k_h \geq 2$. Write $N_i$ (resp., $N_h$) for the level of $f_i$ (resp., $h$), and denote by $\phi_i$ (resp., $\eta$) the nebentypues of $f_i$ (resp., $h$). 
We shall always assume that
\[ k_1 \geq 2-t + k_2+\cdots +k_t + k_h. \]

Let $K$ be a number field obtained by adjoining all the Fourier coefficients of the $f_i$'s and $h$, and the values of the $\phi_i$'s and $\eta$. Fix a prime $\mathfrak{p}$ of $K$ above $p$. For each $\mathfrak{f}\in\{f_1,...,f_t, h\}$, we let $V_{\mathfrak{f}}$ denote the corresponding two dimensional $K_{\mathfrak{p}}$-linear Galois representation attached to $\mathfrak{f}$ in the sense of Deligne.
Writing $\Op$ for the ring of integers of $K_{\mathfrak{p}}$, we denote by $T_{\mathfrak{f}}$ the $\Gal(\overline{\Q}/\Q)$-stable $\Op$-lattice in $V_{\mathfrak{p}}$ defined as in \cite[Section 8.3]{K}. We then set $A_{\mathfrak{f}} = V_{\mathfrak{f}}/T_{\mathfrak{f}}$.

Fix a uniformizer $\varpi$ of $\Op$ and write $k=\Op/\varpi$. The following assumption will be in full force for the remainder of the section.

\begin{enumerate}
 \item[(\textbf{Eis})] There is a short exact sequence of $\Gal(\overline{\Q}/\Q)$-modules
\[ 0\lra \overline{C}_1 \lra A_h[\varpi]\lra \overline{C}_2 \lra  0, \]
where $\overline{C}_1$ and $\overline{C}_2$ are (free) $k$-modules of dimension one.
\end{enumerate}

Let $\overline{\xi}_i: \Gal(\overline{\Q}/\Q) \lra k^\times$ be the character attached to $\overline{C}_i$. Following \cite[Lemma 2.1]{JSV}, we fix a choice of a Dirichlet character $\xi_i$ such that $\xi_i$ is a lift of $\overline{\xi}_i$. Let $C_i$ be the module which is $K/\Op$ as an $\Op$-module with an action of $\Gal(\overline{\Q}/\Q)$ via $\xi_i$. Note that $C_i[\varpi]= \overline{C}_i$.

 Since $f_1$ is assumed to be $p$-ordinary, there exists a one-dimensional $K_{\mathfrak{p}}$-subspace $V_{f_1,p}$ of $V_{f_1}$ invariant under the action of $\Gal(\overline{\Q}_p/\Qp)$ (for instance, see \cite[Proposition 17.1]{K}). Set $A_{f_1,p} = V_{f_1,p}/(T_{f_1}\cap V_{f_1,p})$.

Following \cite{Del, JSV}, we consider the following three sets of data.

\begin{enumerate}
  \item[$(I)$] $B = A_{f_1}\ot_{\Op}\cdots \ot_{\Op} A_{f_t}\ot_{\Op} A_h$, $B_p = A_{f_1,p}\ot_{\Op}\cdots \ot_{\Op} A_{f_t}\ot_{\Op} A_h$;
  \item[$(II)$] $B = A_{f_1}\ot_{\Op}\cdots \ot_{\Op} A_{f_t}\ot_{\Op} C_1$, $B_p = A_{f_1,p}\ot_{\Op}\cdots \ot_{\Op} A_{f_t}\ot_{\Op} C_1$;
  \item[$(III)$] $B = A_{f_1}\ot_{\Op}\cdots \ot_{\Op} A_{f_t}\ot_{\Op} C_2$, $B_p = A_{f_1,p}\ot_{\Op}\cdots \ot_{\Op} A_{f_t}\ot_{\Op} C_2$.
\end{enumerate}

Let $F_\infty$ be an admissible $p$-adic Lie extension of $\Q$. Let $S$ be a finite set of primes of $\Q$ which contains the prime $p$, the prime divisors of $N_1\cdots N_tN_h$, the ramified primes of $F_\infty/F$ and the infinite prime. Denote by $S_1$ the set of primes of $S$ not divisible by $p$ and whose decomposition group in $G$ is of dimension one.

We then set $\Sel_{Gr}(A_{\underline{f}\ot h}/F_\infty)$, $\Sel_{Gr}(A_{\underline{f}(\xi_1)}/F_\infty)$ and $\Sel_{Gr}(A_{\underline{f}(\xi_2)}/F_\infty)$ for the Greenberg Selmer groups of the respective three data.
Similarly, we use $\Sel_{str}(-/F_\infty)$ to denote the respective strict Selmer groups.

For $B\in \{A_{\underline{f}\ot h}, A_{\underline{f}(\xi_1)}, A_{\underline{f}(\xi_2)}\}$, if  $X_{str}(B/F_\infty)\in \N_H(G)$, we set
\[  \tau(B) := \pi_1\big(\ch_{G,\ga}\big(X_{str}(B/F_\infty)\big)\big) \pi_1\big(\ch_{G,\ga}\big(B(F_\infty)^\vee\big)\big)^{-1} \pi_1\left(\ch_{G,\ga}\big(\mathrm{Coind}^{G}_{G_p}\big(B^-_p(F_{\infty,p})\big)^\vee\big)\right)\]
\[\hspace{1.5in} 
\times\prod_{v\in S_1}\pi_1\left(\ch_{G,\ga}\big(\mathrm{Coind}^{G}_{G_v}\big(B(F_{\infty,v})^\vee\big)\big)\right) \]

\[ \hspace{.5in} = \pi_1\big(\ch_{G,\ga}\big(X_{Gr}(B/F_\infty)\big)\big) \pi_1\big(\ch_{G,\ga}\big(B(F_\infty)^\vee\big)\big)^{-1} \pi_1\left(\ch_{G,\ga}\big(\mathrm{Coind}^{G}_{G_p}\big(B^-_p(F_{\infty,p})\big)^\vee\big)\right)\]
\[\hspace{1.5in} 
\times~~\pi_1\left(\ch_{G,\ga}\big(\mathrm{Coind}^{G}_{G_p}H^1(\Gal(F_{\infty,p}^{ur}/F_{\infty,p}),  B_p^-(F_{\infty,p}^{ur}))^\vee\big)\right) \]
\[\hspace{1.5in}\times \prod_{v\in S_1}\pi_1\left(\ch_{G,\ga}\big(\mathrm{Coind}^{G}_{G_v}\big(B(F_{\infty,v})\big)^\vee\big)\right)\quad\in K_1(k\ps{G}_{\overline{\Si}}).\]

We can now state our main theorem.

\bt \label{main}
Suppose that there exist a finite extension $L$ of $\Q$ contained in $F_\infty$ such that $\Gal(F_\infty/L)$ is pro-$p$ and that $X_{str}(A_{\underline{f}(\xi_1)}/L_\cyc)$ and  $X_{str}(A_{\underline{f}(\xi_2)}/L_\cyc)$ are finitely generated over $\Op$. Then $X_{str}(A_{\underline{f}(\xi_1)}/F_\infty)$, $X_{str}(A_{\underline{f}(\xi_2)}/F_\infty)$ and $X_{str}(A_{\underline{f}\ot h}/F_\infty)$ are finitely generated over $\Op\ps{H}$. Furthermore, we have
\[ \tau(A_{\underline{f}\ot h}) = \tau(A_{\underline{f}(\xi_1)}) \tau(A_{\underline{f}(\xi_2)})  \]
in $K_1(k\ps{G}_{\overline{\Si}})$.
\et

\bc \label{maincor}
Retain the assumptions of Theorem \ref{main}. Suppose the following hold.
\begin{enumerate}
  \item[$(a)$] For each prime $v\in S$, the decomposition group $G_v$ of $G$ has dimension at least 2.
  \item[$(b)$] Assume that either one of the following statement is valid.
      \begin{enumerate}
        \item[$(i)$] $A_{\underline{f}(\xi_1)}(F_{\infty})$, $A_{\underline{f}(\xi_2)}(F_{\infty})$, $A_{\underline{f}\ot h}(F_{\infty})$, $A^-_{\underline{f}(\xi_1)}(F_{\infty,p})$ and $A^-_{\underline{f}(\xi_2)}(F_{\infty,p})$ and $A^-_{\underline{f}\ot h}(F_{\infty,p})$ are finite.
        \item[$(ii)$] $G_S(F_\infty)$ acts trivially on  $A_{\underline{f}(\xi_1)}$, $A_{\underline{f}\ot h}$ and $A_{\underline{f}(\xi_1)}$.
      \end{enumerate}
                                             \end{enumerate}
Then we have
\[ \pi_1\big(\ch_{G,\ga}\big(X_{Gr}(A_{\underline{f}\ot h}/F_\infty)\big)\big) = \pi_1\big(\ch_{G,\ga}\big(X_{Gr}(A_{\underline{f}(\xi_1)}/F_\infty)\big)\big) \pi_1\big(\ch_{G,\ga}\big(X_{Gr}(A_{\underline{f}(\xi_2)}/F_\infty)\big)\big).  \]
Furthermore, every open normal pro-$p$ subgroup $H'$ of $H$, we have
\[\rank_{\Op\ps{H'}}\big(X_{Gr}(A_{\underline{f}\ot h}/F_\infty)\big)\big) = \rank_{\Op\ps{H'}}\big(X_{Gr}(A_{\underline{f}(\xi_1)}/F_\infty)\big)\big) + \rank_{\Op\ps{H'}}\big(X_{Gr}(A_{\underline{f}(\xi_2)}/F_\infty)\big)\big).  \]
\ec

\br
It is also possible to compare the $\Op\ps{H'}$-ranks of the Selmer groups in the general context of Theorem 4.1. In this case, the formula will also consist of local terms, and be less pleasant-looking. 
\er 

To prepare for the proof, it is convenient to introduce the following notion.

\bd
Let $\varphi :M\lra N$ be a homomorphism of $k\ps{G}$-modules such that $\ker \varphi$ and $\coker\varphi$ are finitely generated over $k\ps{H}$. We define the content of $\varphi$ to be
\[ c(\varphi) = [\ker \varphi] - [\coker \varphi] \in K_0(\overline{\N}_H(G)).\]

\ed

\bl \label{content}
The following statements are valid.
\begin{enumerate}
  \item[$(a)$] If $\varphi :M\lra N$ is a homomorphism of $k\ps{G}$-modules with $M$ and $N$ being finitely generated over $k\ps{H}$. Then one has
\[ c(\varphi) = [M]- [N]. \]
  \item[$(b)$] Suppose that we are given a commutative diagram\[ \entrymodifiers={!! <0pt, .8ex>+} \SelectTips{eu}{}\xymatrix{
     0 \ar[r] &  M_1\ar[r] \ar[d]_{\varphi_1} &  M_2 \ar[r] \ar[d]_{\varphi_2}&  M_3\ar[d]_{\varphi_3}\ar[r] & 0\\
     0 \ar[r] & N_1\ar[r]_{} &  N_2 \ar[r] &  N_3\ar[r] & 0
     }\]
of $k\ps{G}$-modules with exact rows. If two of $c(\varphi_1)$, $c(\varphi_2)$ and $c(\varphi_3)$ are well-defined, then so is the remaining one and we have
\[ c(\varphi_2)= c(\varphi_1)+ c(\varphi_3). \]
          \end{enumerate}
\el

\bpf
This is a straightforward exercise which is left to the readers.
\epf

We can now give the proof of Theorem \ref{main}.

\bpf[Proof of Theorem \ref{main}]
We first note that it follows from Lemma \ref{fg NHG} that $X_{str}(A_{\underline{f}(\xi_1)}/F_\infty)$ and  $X_{str}(A_{\underline{f}(\xi_2)}/F_\infty)$ lie in $\N_H(G)$. 
To simplify notation, we write $\bar{A}_1 = A_{\underline{f}(\xi_1)}[\varpi], \bar{A}_2 = A_{\underline{f}(\xi_2)}[\varpi]$, $\bar{A}_0 = A_{\underline{f}\ot h}[\varpi]$ and $C_p(i)= \mathrm{Coind}^G_{G_p}(\bar{A}^-_i(F_{\infty,p}))$ and $C_{p'}(i)= \displaystyle\bigoplus_{v\in S_1}\mathrm{Coind}^G_{G_v}(\bar{A}_i(F_{\infty,v}))$. Consider the following commutative diagram

\[
\entrymodifiers={!! <0pt, .8ex>+} \SelectTips{eu}{}\xymatrix{
     0 \ar[r] & \bar{A}_1(F_\infty) \ar[r] \ar[d]_{\be_1} & \bar{A}_0(F_\infty) \ar[r] \ar[d]_{\be_0}& \bar{A}_2(F_\infty) \ar[d]_{\be_2} \ar[r] &&\\
     0 \ar[r] &  C_p(1) \times C_{p'}(1)\ar[r]_{} &   C_p(0) \times C_{p'}(0)  \ar[r] &  C_p(2) \times C_{p'}(2)  \ar[r] & &
     }
\]

\[ \entrymodifiers={!! <0pt, .8ex>+} \SelectTips{eu}{}\xymatrix{
&\ar[r] & H^1(G_S(F_\infty), \bar{A}_1) \ar[d]_{\la_1}  \ar[r] & H^1(G_S(F_\infty), \bar{A}_0) \ar[r] \ar[d]_{\la_0} &  H^1(G_S(F_\infty), \bar{A}_2) \ar[r] \ar[d]_{\la_2} &  0\\
&\ar[r] & \displaystyle\bigoplus_{v\in S_1}J_v(\bar{A}_1/F_\infty) \ar[r]& \displaystyle\bigoplus_{v\in S_1}J_v(\bar{A}_0/F_\infty) \ar[r] & \displaystyle\bigoplus_{v\in S_1}J_v(\bar{A}_2/F_\infty) \ar[r] & 0 
     }
\]
where the top and bottom rows come from the long exact sequence of
\[ 0\lra \bar{A}_1\lra \bar{A}_0\lra \bar{A}_2\lra 0, \] 
and the vertical maps are the usual global-to-local maps. Note that the rightmost zero in the top row is a consequence of Proposition \ref{p-adic torsion}(a) and the rightmost zero in the bottom row is a consequence of \cite[Theorem 7.1.8(i)]{NSW}. Also, one has $\ker \la_i =\Sel_{str}(\bar{A}_i/F_\infty)$ by Proposition \ref{p-adic torsion}(b). Finally, we note that the leftmost three vertical maps are injective homomorphisms of cofinitely generated $k\ps{H}$-modules. By hypothesis, the content of $\la_1$ and $\la_2$ are well-defined. Hence we may apply Lemma \ref{content}(b) to see that the content of $\la_0$ is well-defined and that $\Sel_{str}(\bar{A}_i/F_\infty)$ is cofinitely generated over $k\ps{H}$. In particular, the latter implies that $X_{str}(A_{\underline{f}\ot h}/F_\infty)$ is finitely generated over $\Op\ps{H}$. Furthermore, we have the following equality 
\[ c(\be_1^\vee)  +c(\be_2^\vee) + c(\la_0^\vee) =  c(\be_0^\vee)+ c(\la_1^\vee) + c(\la_2^\vee).  \]
As each $\be_i$ is a homomorphism of finitely generated $k\ps{H}$-modules, we have 
\[ c(\be_i^\vee) =  [C_p(i)^\vee]+[C_{p'}(i)^\vee]- [\bar{A}_i(F_\infty)^\vee]  \]
by Lemma \ref{content}(a).
Combining all these observations with Proposition \ref{char mod p}, followed by a direct calculation, we obtain the required conclusion.
\epf

It remains to prove Corollary \ref{maincor}.

\bpf[Proof of Corollary \ref{maincor}]
Since (a) holds, we see that $S_1=\emptyset$ and so there is no contribution from these terms. We split our argument into two cases.

\underline{Suppose that (a) and (b)(i) hold}

We shall show that under these hypotheses,  one has
\[\tau(B) =\pi_1\big(\ch_{G}\big(X_{Gr}(B/F_\infty)\big)\big)\]
for each $B\in \{A_{\underline{f}\ot h}, A_{\underline{f}(\xi_1)}, A_{\underline{f}(\xi_2)}\}$.  
For the verification of which, it suffices to show that the other terms appearing in the definition of $\tau$ are trivial. Since $B(F_\infty)$ is finite, we may apply \cite[Lemma 2.2]{BZ} to conclude that $[B(F_\infty)^\vee]=0$ in $\N_H(G)$. Hence $\pi_1\big(\ch_{G}\big(B(F_\infty)^\vee\big)\big)$ is trivial in $K_1(k\ps{G}_{\overline{\Sigma}})$. 
Set $H_p := H\cap G_p$. Then one checks easily that $G_p/H_p$ is a nontrivial subgroup of $G/H$ and so a similar argument as above shows that $B^-_p(F_{\infty,p})^\vee$ has trivial class in $\N_{H_p}(G_p)$. It then follows from this that $\mathrm{Coind}^{G}_{G_p}\big(B^-_p(F_{\infty,p})\big)^\vee$ has trivial class in $\N_{H}(G)$ which in turn implies that $\pi_1\left(\ch_{G}\big(\mathrm{Coind}^{G}_{G_p}\big(B^-_p(F_{\infty,p})\big)^\vee\big)\right)$ is trivial in $K_1(k\ps{G}_{\overline{\Sigma}})$. Finally, it remains to show that \[\pi_1\left(\ch_ {G}\big(\mathrm{Coind}^{G}_{G_p}H^1(\Gal(F_{\infty,p}^{ur}/F_{\infty,p}),  B_p^-(F_{\infty,p}^{ur}))^\vee\big)\right)\] is trivial in $K_1(k\ps{G}_{\overline{\Sigma}})$. For this, it again suffices to show that $H^1\big(\Gal(F_{\infty,p}^{ur}/F_{\infty,p}),  B_p^-(F_{\infty,p}^{ur})\big)$ is finite. Let $K_\infty$ be the maximal pro-$p$ extension of $F_{\infty,p}$ contained in $F_{\infty,p}^{ur}$, and write $U = \Gal(K_\infty/F_{\infty,p})$. Since $\Gal(F_{\infty,p}^{ur}/K_\infty)$ is coprime to $p$, one has
\[ H^1(\Gal(F_{\infty,p}^{ur}/F_{\infty,p}),  B_p^-(F_{\infty,p}^{ur})) \cong H^1(U,  B_p^-(K_\infty)). \]

 On the other hand, it follows from \cite[Proposition 5.3.20]{NSW} that
\[ \rank_{\Op\ps{U}} B_p^-(K_\infty)^\vee = \rank_{\Op} B_p^-(F_{\infty,p})^\vee- \rank_{\Op} H^1(U,  B_p^-(K_\infty))^\vee. \]
Since $B_p^-(F_{\infty,p})$ is finite and $B_p^-(K_\infty)^\vee$ is torsion over $\Op\ps{U}$, we have that $H^1(U,  B_p^-(F_{\infty,p}^{ur}))$ is finite as required.

\medskip
\underline{Suppose that (a) and (b)(ii) hold}

In view of (b)(ii), we have $B(F_\infty) = B$ and $B_p^-(F_{\infty,p}) =  B_p^-(F_{\infty,p}^{ur})=B_p^-$ for  $B\in \{A_{\underline{f}\ot h}, A_{\underline{f}(\xi_1)}, A_{\underline{f}(\xi_2)}\}$. Since $B$ is $\Op$-divisible, we have $B^\vee[\varpi]=0$ and so it follows from Corollary \ref{chicommutecoro} that
\[\pi_1\big(\ch_{G}\big(B^\vee\big)\big) = \overline{\ch}_{G}\big(B[\varpi]^\vee\big). \]
 On the other hand, it follows from the natural short exact sequence
\[ 0 \lra  A_{\underline{f}(\xi_1)}[\varpi] \lra A_{\underline{f}\ot h}[\varpi] \lra  A_{\underline{f}(\xi_2)} [\varpi] \lra 0 \]
that 
\[ \overline{\ch}_{G}\big(A_{\underline{f}\ot h}[\varpi]^\vee\big) = \overline{\ch}_{G}\big(A_{\underline{f}(\xi_1)}[\varpi]^\vee\big)
\overline{\ch}_{G}\big(A_{\underline{f}(\xi_2)}[\varpi]^\vee\big). \]
A similar conclusion can be derived for the local terms associated to $B_p^-$. Hence these terms cancel off in the expression of $\tau$ in Theorem \ref{main} leaving behind the Greenberg Selmer groups' characteristics elements which is the asserted equality of the corollary.
\epf


\footnotesize
\begin{thebibliography}{00}

\bibitem{AAS} S. Ahmed, C. Aribam and S. Shekhar, Root numbers and parity of local Iwasawa invariants, J. Number Theory 177 (2017), 285–306.
\bibitem{BZ} T. Backhausz and G. Z\'abr\'adi, Algebraic functional equations and
completely faithful Selmer groups, Int. J. Number Theory
11 (2015) 1233-1257.

\bibitem{BerK} A. J. Berrick and M. E. Keating, The localization sequence in $K$-theory, J. K-Theory 9 (1995), no. 6, 577-589.

\bibitem{Bu09} D. Burns, Algebraic $p$-adic $L$-functions in non-commutative Iwasawa theory. Publ. Res. Inst. Math. Sci. 45 (2009), no. 1, 75-87.

\bibitem{Bu15} D. Burns, On main conjectures in non-commutative Iwasawa theory and related conjectures. J. Reine Angew. Math. 698 (2015), 105-159.

\bibitem{BV} D. Burns and O. Venjakob, On descent theory and main conjectures in non-commutative Iwasawa theory. J. Inst. Math. Jussieu 10 (2011), no. 1, 59-118.

\bibitem{CDLSS} J. Coates, T. Dokchitser, Z. Liang, W. Stein and R. Sujatha,
Non-commutative Iwasawa theory for modular forms, Proc. London Math. Soc. (3) 107 (2013) 481-516.


\bibitem{CF+} J. Coates, T. Fukaya, K. Kato, R. Sujatha and O. Venjakob, The $GL_2$ main conjecture for elliptic curves without complex multiplication.
Publ. Math. Inst. Hautes \'Etudes Sci. No. 101 (2005), 163-208.




\bibitem{CS12} J. Coates and R. Sujatha, On the $\M_H(G)$-conjecture, in: Non-abelian fundamental groups and Iwasawa theory, 132-161, London Math. Soc. Lecture Note Ser., 393, Cambridge Univ. Press, Cambridge, 2012. 

\bibitem{DelPLMS} D. Delbourgo,
Variation of the analytic $\la$-invariant over a solvable extension, Proc. London Math. Soc. (3) 120 (2020) 918-960.

\bibitem{Del} D. Delbourgo, On the Iwasawa $\mu$-invariant and $\lambda$-invariant associated to tensor products of newforms, Ann. Inst. Fourier (Grenoble) 74 (2024), no. 2, 451-502.

\bibitem{DelGil} D. Delbourgo and H. Gilmore, Controlling $\lambda$-invariants for the double and triple product $p$-adic
$L$-functions, J. Th\'eor. Nombres Bordeaux  33, no 3.1 (2021), 733-778. 

\bibitem{EPW} M. Emerton, R. Pollack and T. Weston, Variation of Iwasawa invariants in Hida families, Invent. Math., 163 (2006) 523-580. 
    
\bibitem{FK} T. Fukaya and K. Kato, A formulation of conjectures on $p$-adic zeta functions in noncommutative Iwasawa theory. Proceedings of the St. Petersburg Mathematical Society. Vol. XII, 1-85, Amer. Math. Soc. Transl. Ser. 2, 219, Amer. Math. Soc., Providence, RI, 2006.

\bibitem{GW} K. R. Goodearl and R. B. Warfield, An
introduction to non-commutative Noetherian rings, London Math. Soc.
Stud. Texts 61, Cambridge University Press, 2004.

 \bibitem{G89} R. Greenberg, Iwasawa theory for $p$-adic representations, in:
Algebraic Number Theory-in honor of K. Iwasawa, ed. J.
Coates, R. Greenberg, B. Mazur and I. Satake, Adv. Std. in Pure
Math. 17, 1989, pp. 97-137.

\bibitem{GV} R. Greenberg and V. Vatsal, On the Iwasawa invariants of elliptic curves, Invent. Math., 142 (2000) 17–63.

\bibitem{Ha} Y. Hachimori, Iwasawa $\lambda$-invariants and congruence of Galois representations, J. Ramanujan Math. Soc., 26(2) (2011) 203–217.  

\bibitem{How} S. Howson, Euler characteristics as invariants of Iwasawa modules, Proc. Lond. Math. Soc. 85 (3) (2002) 634-658. 

\bibitem{Iw73} K. Iwasawa, On $\Z_l$-extensions of algebraic number fields. Ann. of Math. (2) 98 (1973), 246-326.

\bibitem{jannsen} U. Jannsen, A spectral sequence for Iwasawa adjoints, M\"{u}nster J. Math. 7 (2014), no. 1, 135-148.

\bibitem{JMS} S. Jha, T. Mandal and S. Shekhar, Multiplicities in Selmer groups and root numbers of Artin twists, J. Number Theory  238 (2022), 147-182.

\bibitem{JSV} S. Jha,  S. Shekhar and R. Vangala, Iwasawa theory for Rankin-Selberg product at an Eisenstein prime,
 arXiv:2209.04482 [math.NT].

\bibitem{Kak} M. Kakde, The main conjecture of Iwasawa theory for totally
real fields. Invent. Math, 193 (2013) 539-626.

\bibitem{K} K. Kato, $p$-adic Hodge theory and values of zeta functions of
modular forms, in: Cohomologies $p$-adiques et applications
arithm\'etiques. III., Ast\'erisque 295, 2004, ix, pp.
117-290.

\bibitem{Lam} T. Y. Lam, Lectures on Modules and Rings.
Grad. Texts in Math. 189, Springer, 1999.


\bibitem{LimFine}  M. F. Lim, Notes on the fine Selmer groups. Asian J. Math. 21 (2017), no. 2, 337-362.

\bibitem{LimCMu} M. F. Lim, Comparing the $\pi$-primary
   submodules of the dual Selmer groups.
 Asian J. Math. 21 (2017), no. 6, 1153-1182.
 
 
\bibitem{Maz} B. Mazur, Rational points of abelian varieties with values in towers of number fields, Invent. Math. 18 (1972) 183-266.

\bibitem{Mi} J. Milne, Arithmetic Duality Theorems. Second edition. BookSurge, LLC, Charleston, SC, 2006. viii+339 pp.

\bibitem{NSW} J. Neukirch, A. Schmidt and K. Wingberg,
Cohomology of Number Fields, 2nd edn., Grundlehren Math.
Wiss. 323 (Springer-Verlag, Berlin, 2008).

\bibitem{Neu} A. Neumann, Completed group algebras without zero divisors.
Arch. Math. 51(6) (1988) 496-499.

\bibitem{OV} Y. Ochi and O. Venjakob, On the ranks of Iwasawa modules over $p$-adic Lie extensions. Math. Proc. Camb. Phil. Soc. 135 (2003), 25-43.

\bibitem{Sch85} P. Schneider, $p$-adic height pairings II. Invent. Math. 79 (1985), no. 2, 329-374.

\bibitem{Sh} S. Shekhar, Parity of ranks of elliptic curves with equivalent mod $p$ Galois representations, Proc. Amer. Math. Soc. 144 (2016), 3255–3266.
    
\bibitem{V02} O. Venjakob, On the structure theory of the Iwasawa algebra
of a $p$-adic Lie group. J. Eur. Math. Soc. 4(3)
(2002) 271-311.

\bibitem{V03} O. Venjakob, A non-commutative Weierstrass preparation theorem and applications to Iwasawa theory. J. Reine Angew. Math. 559 (2003) 153-191.

\bibitem{V05} O. Venjakob, Characteristic elements in noncommutative Iwasawa theory. J. Reine Angew. Math. 583 (2005), 193-236.


\bibitem{V07} O. Venjakob, From the Birch and Swinnerton-Dyer conjecture to non-commutative Iwasawa theory via the equivariant Tamagawa number conjecture-a survey. $L$-functions and Galois representations, 333-380, London Math. Soc. Lecture Note Ser., 320, Cambridge Univ. Press, Cambridge, 2007

\bibitem{Wei} C. Weibel, The $K$-book. An introduction to algebraic $K$-theory. Graduate Studies in
Mathematics, 145. American Mathematical Society, Providence, RI, 2013.

\bibitem{Wi} M. Witte,  On a localisation sequence for the $K$-theory of skew power series rings, J. K-Theory 11 (2013), no. 1, 125-154.

\end{thebibliography}
\end{document}